\documentclass[10pt]{amsart}

\usepackage{palatino}
\usepackage{mathpazo}

\usepackage{amssymb}
\usepackage{amsmath}
\usepackage{amscd}

\usepackage[alphabetic]{amsrefs}
\usepackage{mathrsfs}

\newcommand{\sep}{{\operatorname{sep}}}
\newcommand{\alg}{{\operatorname{alg}}}
\newcommand{\FFF}{\mathfrak{f}}

\renewcommand{\setminus}{\ \rule[2.5pt]{7pt}{1pt}\ }

\newcommand{\glie}{\mathfrak{g}}
\newcommand{\hlie}{\mathfrak{h}}

\newcommand{\Lie}{\operatorname{Lie}}
\newcommand{\G}{\mathbf{G}}
\newcommand{\Gm}{{\G_m}}

\newcommand{\Aff}{\mathbf{A}}
\newcommand{\Proj}{\mathbf{P}}
\newcommand{\Ad}{\operatorname{Ad}}
\newcommand{\ad}{\operatorname{ad}}
\newcommand{\Int}{\operatorname{Int}}

\newcommand{\Der}{\operatorname{Der}}

\newcommand{\Hom}{\operatorname{Hom}}

\newcommand{\congruent}{\equiv}

\newcommand{\Par}{\underline{\operatorname{Par}}}

\newcommand{\B}{\mathcal{B}}
\newcommand{\A}{\mathcal{A}}
\newcommand{\C}{\mathcal{C}}
\newcommand{\Ah}{\A^{\operatorname{h}}}
\newcommand{\Ash}{\A^{\operatorname{sh}}}

\newcommand{\F}{\mathcal{F}}
\newcommand{\FF}{\mathbf{F}}

\newcommand{\Et}[1]{\operatorname{Et}_{/#1}}
\newcommand{\et}{\operatorname{et}}

\newcommand{\OO}{\mathscr{O}}

\newcommand{\LL}{\mathcal{L}}

\newcommand{\Z}{\mathbf{Z}}
\newcommand{\Q}{\mathbf{Q}}

\newcommand{\pp}{\mathfrak{p}}

\newcommand{\SL}{\operatorname{SL}}
\newcommand{\GL}{\operatorname{GL}}

\newcommand{\Spec}{\operatorname{Spec}}
\newcommand{\tensor}{\otimes}
\newcommand{\Stab}{\operatorname{Stab}}
\newcommand{\Trans}{\operatorname{Trans}}

\swapnumbers

\theoremstyle{plain}

\newtheorem*{theorem}{Theorem}

\swapnumbers
\newtheorem{theoreml}{Theorem}

\swapnumbers

\newtheorem*{prop}{Proposition}
\newtheorem*{cor}{Corollary}

\swapnumbers

\swapnumbers

\swapnumbers
\newtheorem{stmt}{}[subsection]

\swapnumbers
\theoremstyle{remark}
\newtheorem*{rem}{Remark}

\numberwithin{equation}{subsection}

\begin{document}

\title{The centralizer of a nilpotent section}
\author{George J. McNinch}
\address{Department of Mathematics,
         Tufts University,
         503 Boston Avenue,
         Medford, MA 02155,
         USA}
\email{george.mcninch@tufts.edu, mcninchg@member.ams.org}
\date{June 9, 2007}
\thanks{\noindent Research of the author supported
  in part by the US National Science Foundation through DMS-0437482.}

\dedicatory{To Toshiaki Shoji, with respect and admiration, on the
  occasion of his 60th birthday.}

\begin{abstract}
  Let $F$ be an algebraically closed field and let $G$ be a semisimple
  $F$-algebraic group for which the characteristic of $F$ is
  \emph{very good}.  If $X \in \Lie(G) = \Lie(G)(F)$ is a nilpotent element in
  the Lie algebra of $G$, and if $C$ is the centralizer in $G$ of $X$,
  we show that (i) the root datum of a Levi factor of $C$, and (ii)
  the component group $C/C^o$ both depend only on the Bala-Carter
  label of $X$; i.e. both are independent of very good characteristic.
  The result in case (ii) depends on the known case when $G$ is
  (simple and) of adjoint type.

  The proofs are achieved by studying the centralizer $\C$ of a
  nilpotent section $X$ in the Lie algebra of a suitable semisimple
  group scheme over a Noetherian, normal, local ring $\A$. When the
  centralizer of $X$ is equidimensional on $\Spec(\A)$, a crucial
  result is that locally in the \'etale topology there is a smooth
  $\A$-subgroup scheme $L$ of $\C$ such that $L_t$ is a Levi factor of
  $\C_t$ for each $t \in \Spec(\A)$.


\end{abstract}

\maketitle
\setcounter{tocdepth}{1}
\tableofcontents

\section{Introduction}

\subsection{The main results}
\label{sub:main-results-intro}

Let $E$ and $F$ be algebraically closed fields, and let $G_E$ and
$G_F$ be semisimple algebraic groups over $E$ and $F$ respectively.
We are going to assume that the \emph{root data} of these two groups
coincide. Further, we suppose that the characteristic of $E$ is 0, and
that the characteristic of $F$ is \emph{very good} for $G_F$ -- see
\S\ref{sub:good-primes}.

Using the Bala-Carter Theorem \ref{stmt:Bala-Carter-theorem}, we may
identify the set of nilpotent orbits of $G_F$ in $\Lie(G_F) = \glie_F$ with
the set of nilpotent orbits of $G_E$ in $\Lie(G_E) = \glie_E$.

Suppose that the orbits of the nilpotent elements $X_E \in \glie_E$
and $X_F \in \glie_F$ are the same under the Bala-Carter
identification, let $C_E$ be the centralizer of $X_E$ in $G_E$,
and let $C_F$ be the centralizer of $X_F$ in $G_F$.

If $H$ is an algebraic group, one says that a closed subgroup $L
\subset H$ is a Levi factor if the connected component $L^o$ is
reductive and if $H$ is isomorphic as an algebraic group to the
semidirect product $L \cdot R_u(H)$, where we have written $R_u(H)$
for the unipotent radical of $H$.

The groups $C_F$ and $C_E$ have Levi factors $L_F \subset C_F$ and
$L_E \subset C_E$. Indeed, this is immediate in characteristic 0, since
a result of Mostow shows \emph{every} linear group to have a Levi
factor in that case; in positive characteristic, existence of a Levi
factor for $C_F$ may be deduced as a consequence of Premet's recent
conceptual proof of the Bala-Carter theorem \cite{premet}; cf.
\ref{stmt:assoc-cochar} below.  The main results of this paper may now
be stated:

\begin{theoreml}
  \label{theorem:levi-factor}
  The root datum of the reductive group $L_F^o$ may be identified with
  that of $L_E^o$.
\end{theoreml}

\begin{theoreml}
  \label{theorem:component}
  The finite groups $C_E/C_E^o$  and $C_F/C_F^o$ are isomorphic.
\end{theoreml}

When $p$ is not a good prime for $G$, the Bala-Carter parametrization
of nilpotent orbits does not hold; cf. \cite{carter}*{\S5.11} and
\cite{jantzen-nil}*{\S5.13 -- 5.15} for examples of ``extra''
nilpotent orbits for these primes. So our statements must at least
exclude ``bad'' characteristics. We have not, however, attempted to
prove our results for semisimple groups in all good
characteristics. Instead, we have chosen to prove the theorems of this
paper under some ``standard'' assumptions on $G$; in fact, we will
prove Theorems \ref{theorem:levi-factor} and \ref{theorem:component}
for the \emph{$T$-standard} reductive groups introduced in \S
\ref{sub:strongly-standard}. A semisimple group is $T$-standard in
case the characteristic is \emph{very good} for $G$, but the group
$\GL_n$ is always $T$-standard. Thus, our statements apply, for
example, to the group $\GL_n$ for any $n$, but not to $\SL_n$ when $n
\equiv 0 \pmod p$.  Note that the centralizer of a \emph{regular}
nilpotent element in $\SL_n$ is the direct product of a connected
unipotent group with the group $\mu_n$ of $n$-th roots of unity; thus
when $n \equiv 0 \pmod{p}$, the naive statement of Theorem
\ref{theorem:component} would not be correct for $\SL_n$.

In the remainder of this introduction, we give an overview of our
strategy of proofs of Theorems \ref{theorem:levi-factor} and
\ref{theorem:component}.  We first observe that -- as a consequence of
the Bala-Carter Theorem; see \ref{stmt:choose-favorite-fields} -- it
suffices to prove these Theorems after making a particular choice for
the fields $E$ and $F$.  For convenience, we will prove the result
when $F$ is an algebraic closure of the finite field $\mathbf{F}_p$,
and $E$ is some suitable algebraically closed field of characteristic
0.  The proofs will be given in \S \ref{sub:proof-levi-factor} and \S
\ref{sub:proof-component}.  We now give some further details about
these proofs.

\subsection{The instability parabolic}

As already mentioned, we rely on the fact that the nilpotent orbits for the group
$G_E$ and for the group $G_F$ are described by the Bala-Carter
theorem; cf. \S \ref{sub:bala-carter}. 

Recall that a key idea behind Premet's recent proof \cite{premet} of the
Bala-Carter theorem was to use a result in geometric invariant theory
-- due to Kempf and to Rousseau -- which attaches a collection of
optimal cocharacters of $G$ to an unstable vector in a linear
representation of $G$. Let us explain this a bit more. Write $G$
for one of the groups $G_E$ or $G_F$. An element $X \in \Lie(G)$ is
nilpotent if and only if the closure of its adjoint orbit contains 0;
such vectors are said to be unstable. The Hilbert-Mumford criterion for
instability asserts that an unstable vector for $G$ is also unstable
for certain one-dimensional sub-tori of $G$. The Hilbert-Mumford
criterion has a more precise form due to Kempf and to Rousseau: there
is a class of \emph{optimal} cocharacters of $G$ whose images exhibit
such one dimensional sub-tori. One of the nice features of these
optimal cocharacters is that they each define the same parabolic
subgroup $P_X$ of $G$; this parabolic subgroup is known as the
instability parabolic determined by $X$. The instability parabolic
subgroups determined by nilpotent elements play an important role in
this paper.

When $G$ is a reductive group over an arbitrary field $K$ and when $X
\in \Lie(G)(K)$ is nilpotent, one knows e.g. by
\cite{mcninch-rat}*{Prop. 27} that $P_X$ is a $K$-parabolic; cf.
\ref{sub:associated} for more on these matters.

\subsection{The component group when $G$ is of adjoint type}

Let again $G$ be one of the groups $G_E$ or $G_F$, and suppose that
$G$ is of adjoint type.  For a reductive group in characteristic zero,
Alekseevski{\u\i} \cite{alek} determined the structure of the group of
components $C_G(X)/C_G(X)^o$ for each nilpotent $X \in \Lie(G)$.
Sommers \cite{sommers} gave later a more conceptual argument for the
determination of these groups.

Moreover, given nilpotent elements $X_E \in \Lie(G_E)$ and $X_F \in
\Lie(G_F)$ with the same Bala-Carter label, one knows for semisimple
groups of adjoint type that $C_E/C^o_E \simeq C_F/C^o_F$. For a while,
this was known only through case-checking -- especially, by the work
of Mizuno \cite{mizuno}. More recently, this isomorphism was proved by
Premet \cite{premet} and by McNinch-Sommers \cite{mcninch-sommers}.
Thus, the assertion of Theorem \ref{theorem:component} is known
already provided that $G_E$ and $G_F$ are of adjoint type.

\subsection{Group schemes}
The proofs of Theorems \ref{theorem:levi-factor} and
\ref{theorem:component} are achieved by studying reductive group
schemes over more general base schemes. Let us give here a brief
overview of the argument.

We consider a normal, local, Noetherian integral domain $\A$ with
residue field $k$ and field of fractions $K$.  Recall that a point $t
\in \Spec(\A)$ is the same as a prime ideal $\pp_t \subset \A$; write
$k(t)$ for the field of fractions of $\A/\pp_t$. The closed point $s
\in \Spec(\A)$ is the maximal ideal of $\A$, so that $k(s) = k$ is the
residue field. And the generic point $\eta \in \Spec(\A)$ is the prime
ideal $0$ of $\A$, so that $k(\eta) = K$ is the field of fractions.
For any $t \in \Spec(\A)$, we write $k(\bar{t})$ for a separable
closure of the field $k(t)$.

Let $G$ be a semisimple group scheme over $\A$. For $t \in \Spec(\A)$,
we write $G_t$ for the group $G_{/k(t)}$ obtained by base-change; thus
$G_{/k(t)}$ is a semisimple group over the field $k(t)$. We insist
that the characteristic of $k$ is \emph{very good} for $G_s=G_{/k}$;
it is then immediate that the
characteristic of $K$ is very good for $G_\eta = G_{/K}$; see \S
\ref{sub:good-primes}.

An $\A$-section of the Lie algebra $\glie = \Lie(G)$ of $G$ -- i.e.
an element $X \in \glie(\A)$ -- is nilpotent if its image $X_K \in
\glie(K)$ is nilpotent. If $X \in \glie(\A)$ is a nilpotent section,
write $C = C_G(X)$ for the centralizer -- thus $C$ is an $\A$-subgroup
scheme of $G$.  Now, $C_{/K}$ identifies with the centralizer in
$G_{/K}$ of $X_K$, and likewise for $C_{/k}$. If the groups $C_{/K}$
and $C_{/k}$ have the same dimension, we say that $X$ is
equidimensional; we prove in that case -- see Proposition
\ref{sub:smooth-over-A} -- that the group scheme $C$ is smooth over
$\A$. 

If $P_0$ denotes the instability parabolic subgroup of $G_{/K}$
determined by the nilpotent element $X_K \in \glie(K)$, we prove --
see Proposition \ref{sub:instability-parabolic-over-A} -- that there
is a parabolic $\A$-subgroup scheme $P \subset G$ for which $P_{/K} =
P_0$.

This assertion is immediate in case $\A$ is a discrete valuation ring;
see \ref{stmt:parabolic-over-dvr}. The general case is a consequence
of \ref{stmt:section-of-projective-over-normal-A}. Note that for
general $\A$ as above, the conclusion of Proposition
\ref{sub:instability-parabolic-over-A} actually holds by construction
for a collection of equidimensional sections $X \in \glie(\A)$
realized as ``Richardson sections''; see Theorem
\ref{sub:existence-equidim}.

The existence of Richardson sections just mentioned also shows that
for $s \in \Spec(\A)$ and a nilpotent element $Y \in
\glie(k(\overline{s}))$, there is an equidimensional nilpotent section
$X \in \glie(\A)$ such that $Y$ is geometrically conjugate to the
value $X(s) \in \glie(k(s))$ of $X$; moreover, the construction of
this $X$ makes clear that the Bala-Carter datum of the nilpotent
element $X(t) \in \glie(k(\bar{t}))$ is constant for $t \in
\Spec(\A)$.

We may now state a key result: locally in the \'etale topology of
$\Spec(\A)$, the centralizer $C$ has a Levi factor. This means that
after possibly replacing $\A$ by a finite, \'etale, local extension,
we may find a closed, smooth subgroup scheme $L \subset C$ such that
$L^o$ is reductive and such that $L_t$ is a Levi factor of $C_t$
for $t \in \Spec(\A)$; cf. Theorem \ref{sub:max-tori-and-levi}.

The existence of the Levi factor $L$ essentially settles Theorem
\ref{theorem:levi-factor}.  Note that we also prove -- cf.  Corollary
\ref{sub:max-tori-and-levi} -- for \emph{any} equidimensional
nilpotent section $X$ that the Bala-Carter datum of $X(t)$ is constant
for $t \in \Spec(\A)$.

For Theorem \ref{theorem:component} one considers the sheaf on the
\'etale site of $\A$ determined by the quotient $C/C^o \simeq L/L^o$.
When $G$ is a semisimple group scheme over $\A$ with adjoint root
datum, Theorem \ref{theorem:component} is known for the geometric
fibers of $G$; it follows that the sheaf $L/L^o$ is represented by a
finite \'etale group scheme over $\A$. To complete the proof of
Theorem \ref{theorem:component}, we must argue when $G$ is no longer
adjoint that the sheaf $L/L^o$ is still represented by a finite
\'etale group scheme; this is carried out in \S
\ref{section:nilpotent-comp}.


Theorem A was announced by the author in June 2005 in a talk in the
conference on Algebraic Groups and Finite Reductive Groups at the
Bernoulli Center of the \'Ecole Polytechnique F\'ed\'erale de Lausanne. The
author thanks Jens Carsten Jantzen, Michel Raynaud, and Jean-Pierre
Serre for useful remarks during the preparation of this manuscript.

\section{Some recollections}

\subsection{Assumptions and notation}
Let $\A$ be a Noetherian integral domain. We are going to consider
schemes over $\A$, and -- as e.g. in \cite{JRAG} -- we will
interchangeably regard a scheme over $\A$ either as a set-valued
functor on all commutative $\A$-algebras [more precisely: all
commutative $\A$-algebras in some universe, to avoid well-known
logical pitfalls] or as the ringed topological space which represents
this functor.

Given a scheme $X$ of finite type over $\A$ and a point $x \in X$ we
write $\OO_x$ for the local ring of $x$, and we write $k(x)$ for the
residue field of $\OO_x$. When $X =\Spec(\A)$, we write $\A_x$ for
this local ring. We will denote by $k(\bar{x})$ a separably closed
field containing $k(x)$; thus $\bar{x}$ is a \emph{geometric point} of
$X$.

If $t$ is a point of $\Spec(\A)$, we write $X_t$ for the fiber
product $X \times_{\Spec(\A)} \Spec k(t)$; then $X_t$ is a scheme of finite
type over the field $k(t)$.

Similarly, if $\A \subset \B$ is an extension, we write $X_{/\B}$ for the fiber
product $X \times_{\Spec(\A)} \Spec \B$; then $X_{/\B}$ is a scheme
of finite type over $\B$.

\subsection{Normal local rings and \'etale extensions}
\label{sec:normality}

Assume  that the Noetherian
integral domain $\A$ is \emph{normal} and \emph{local}. If $\B$ is a
commutative ring containing $\A$, then $\B$ is said to be a
\emph{local extension} of $\A$ if $\B$ is itself local, and if the
maximal ideal of $\A$ is contained in that of $\B$.  
We have the following:

\begin{stmt} \cite{SGA1}*{Exp. I, Prop. 10.1}
\label{stmt:etale-extension-of-normal}
If $\B \supset \A$ is a local \'etale extension of finite type, then
$\B$ is a domain, the field of fractions $L$ of $\B$ is a finite
separable extension of the field of fractions $K$ of $\A$, and $\B$ is
the integral closure of $\A$ in $L$. In particular, $\B$ is Noetherian
and is finite over $\A$.
\end{stmt}

In the language of the \'etale topology (see \cite{milne}*{II} or \S
\ref{sub:sheaves} below), we have:
\begin{stmt}
  Any \'etale neighborhood $X \to \Spec(\A)$ of the closed point of
  $\Spec(\A)$ contains an affine open $\A$-subscheme $\Spec(\B)
  \subset X$ for some finite, \'etale, local $\A$-algebra $\B$.
\end{stmt}

If $X$ is an $\A$-scheme, then we say that a property of $X$ holds
locally in the \'etale topology of $\A$ if the property holds for the
$\B$-scheme $X_{/\B}$ for a suitable finite, \'etale, local extension
$\B$ of $\A$; note that $\B$ is then necessarily a domain.

\subsection{Smoothness of stabilizers}
\label{sub:smoothness}

In this section, $\A$ is a Noetherian integral domain.  Let $G$ be a
group scheme which is smooth and of finite type over $\A$, and let $Y$
be a scheme which is flat and of finite type over $\A$. Suppose that
$G$ acts on $Y$ and that the action is given by a morphism of
$\A$-schemes $\mathbf{a}:G \times_{\A} Y \to Y$.

If $\alpha \in Y(\A)$ is an $\A$-section, then for each commutative
$\A$-algebra $\Lambda$, the section $\alpha$ determines a section
$\alpha_\Lambda \in Y(\Lambda)$; for $t \in \Spec(\A)$ we write
$\alpha(t)$ for the image $\alpha_{k(t)}$ in $Y(k(t))$.

Let now $\alpha,\beta \in Y(\A)$ be two $\A$-sections of $Y$.  The
transporter $\Trans_G(\alpha,\beta)$ is the subfunctor of $G$ given for each
commutative $\A$-algebra $\Lambda$ by
\begin{equation*}
  \Trans_G(\alpha,\beta)(\Lambda) = \{g \in G(\Lambda) \mid g.\alpha_\Lambda = \beta_\Lambda\}.
\end{equation*}
In particular, the stabilizer $\Stab_G(\alpha) = \Trans_G(\alpha,\alpha)$ is the
subfunctor of $G$ given by
\begin{equation*}
  \Stab_G(\alpha)(\Lambda) = \{g \in G(\Lambda) \mid g.\alpha_\Lambda = \alpha_\Lambda\}.
\end{equation*}

Write $\mu_\alpha$ for the orbit mapping
\begin{equation*}
  \mu_\alpha:G  \to Y
\end{equation*}
determined by the section $\alpha \in Y(\A)$; for each commutative
$\A$-algebra $\Lambda$ we have $\mu_\alpha(h) =
\mathbf{a}(h,\alpha_\Lambda)$ for each $h \in G(\Lambda)$.  We now
regard $\alpha,\beta \in Y(\A)$ as sections $S
\xrightarrow{\phi_\alpha} Y$ and $S \xrightarrow{\phi_\beta} Y$, where
$S$ is the spectrum of $\A$.  Then we see that the sub-functor
$\Trans_G(\alpha,\beta)$ may be identified with the fiber product $G
\times_{Y,\phi_\beta} S$:
\begin{equation*}
  \begin{CD}
    G @<<< G \times_{Y,\phi_\beta} S = T \\
    @V{\mu_\alpha}VV @VVV \\
    Y @<\phi_\beta<< Y \times_{Y,\phi_\beta} S = S.
  \end{CD}
\end{equation*}
It is thus a subscheme of $G$ which is of finite type over $\A$, and
it is closed in $G$ if $\phi_\beta$ is a closed embedding.  In particular,
the stabilizer $\Stab_G(\alpha)$ is a subscheme of $G$ which is of finite
type over $\A$; it is closed in $G$ in case $\phi_\alpha$ is a closed
embedding.

We are interested in conditions under which the transporter
$\Trans_G(\alpha,\beta)$ and the stabilizer $\Stab_G(\alpha)$ are
smooth; we give two such conditions, as follows:

\begin{stmt}
  \label{stmt:smoothness-I}
  Suppose for each $s \in \Spec(\A)$ that the $G_s$-orbit of
  $\alpha(s)$ in $Y_s$ is separable and dense.  Then for each
  section $\beta \in Y(\A)$, the transporter
  $T=\Trans_G(\alpha,\beta)$ is a smooth $\A$-subscheme of $G$. In
  particular, the stabilizer $C=\Stab_G(\alpha)$ is a smooth
  $\A$-subgroup scheme of $G$.
\end{stmt}

\begin{proof}
  If $\mu_\alpha$ denotes the orbit map as above, we claim first that
  $\mu_\alpha:G \to Y$ is smooth.  Fix a point $g \in G$; we argue
  that $\mu_\alpha$ is smooth at $g$.  Indeed, $G$ is smooth over
  $\A$, $Y$ is flat over $\A$, and $\mathbf{a}$ is of finite type.
  Thus according to \cite{SGA1}*{Exp. II, Cor. 2.2}, the
  smoothness of $\mu_\alpha$ at $g$ will follow provided that
  $\mu_{\alpha,s}:G_s \to Y_s$ is smooth at $g$, where
  $s \in \Spec(\A)$ is the image of $g$ under the projection $G \to
  \Spec(\A)$; since the $G_s$-orbit of $\alpha(s)$ is separable
  and dense, $\mu_{\alpha,s}$ is indeed smooth at $g$.

  Since smooth morphisms are stable under base-change \cite{SGA1}*{Exp. II,
    Prop. 1.3}, and since $\mu_\alpha$ is smooth, it follows that $T$
  and $C$ are smooth over $\Spec(\A)$, as required.
\end{proof}

\begin{stmt}
  \label{stmt:smoothness-II}
  Assume that the Noetherian integral domain $\A$ is normal. 
  \begin{enumerate}
  \item If $\Trans_G(\alpha,\beta)_t$ is smooth over $k(t)$ for
    each $t \in \Spec(\A)$, and if the irreducible components of
    $\Trans_G(\alpha,\beta)_t$ all have the same dimension
    independent of $t \in \Spec(\A)$, then $\Trans_G(\alpha,\beta)$ is
    a smooth $\A$-subscheme of $G$.
  \item If $C_G(X)_t$ is smooth over $k(t)$ for each $t \in
    \Spec(\A)$, and if the irreducible components of $C_G(X)_t$
    all have the same dimension independent of $t \in\Spec(\A)$, then
    $C_G(X)$ is a smooth $\A$-subgroup scheme of $G$.
  \end{enumerate}
\end{stmt}

\begin{proof}
  Of course, (2) is a special case of (1).  Since the scheme
  $\Trans_G(X,Z)$ is of finite type over $\A$ and since $\A$ is
  normal, (1) follows from \cite{SGA1}*{Exp. II, Prop. 2.3}.
\end{proof}

\subsection{Henselian rings}
\label{sub:hensel}

Suppose that the Noetherian domain $\A$ is moreover local. Then $\A$
is said to be \emph{Henselian} if the conclusion of Hensel's lemma
holds for $\A$ (see e.g. \cite{milne}*{I \S4}).
The local domain $\A$ is said to be \emph{strictly Henselian} if it is
Henselian and if its residue field is separably closed.  

If $\mathfrak{m}$ is the maximal ideal of $\A$, recall that $\A$ is
Henselian if it is complete in its $\mathfrak{m}$-adic topology
\cite{milne}*{I, Prop. 4.5}.  Given any local domain $\A$, we can
construct its Henselization $\Ah$ \cite{milne}*{I \S4} and its
strict Henselization $\Ash$ \emph{loc. cit.}.

Let $X$ be a smooth scheme of finite type over $\A$, and write $k$ for
the residue field of $\A$.
\begin{stmt}\cite{milne}*{I Exerc. 4.13}
  \label{stmt:henselian-lift-to-A}
  If $\A$ is Henselian, the natural map $X(\A) \to X(k)$ is
  surjective.
\end{stmt}

Suppose that $G$ is a smooth group scheme of finite type over $\A$.
Let $Y$ be a scheme which is flat and of finite type over $\A$, and
assume that $G$ acts on $Y$ by $\A$-morphisms.
\begin{stmt}
  \label{stmt:transport-via-etale-ext}
  Let $\alpha,\beta \in Y(\A)$, suppose that the $G_t$-orbit of
  $\alpha(t)$ is separable and dense in $Y_t$ for every $t \in
  \Spec(\A)$, and suppose that the elements $\alpha(s),\beta(s) \in
  Y(k(s)) = Y(k)$ are conjugate by an element of $G(k)$, where $s$ is
  the closed point of $\Spec(\A)$.  Then locally in the \'etale
  topology of $\A$, the sections $\alpha$ and $\beta$ are conjugate by
  a section of $G$. More precisely, there is a finite, \'etale, local
  extension $\B$ of $\A$ such that $\alpha$ and $\beta$ are conjugate
  by an element of $G(\B)$.
\end{stmt}

\begin{proof}
  Indeed, by \ref{stmt:smoothness-I}, the transporter
  $\Trans_G(\alpha,\beta)$ is a smooth subscheme of $G$. Thus by
  \ref{stmt:henselian-lift-to-A} any $k$-point of the transporter may
  be lifted to an $\Ah$-point, where $\Ah$ is the Henselization. Now
  the result follows from the construction of $\Ah$ as the limit of
  \'etale neighborhoods of $\A$ \cite{milne}*{I \S4}.
\end{proof}


\subsection{Dimensions of fibers}
Let $X$ be a scheme of finite type over the Noetherian domain
$\A$.  Let us write $\pi:X \to S = \Spec(\A)$ for the structure
morphism. Then  Chevalley's upper semicontinuity theorem --
cf.  \cite{EGA-IV}*{Thm. 13.1.3} -- gives rough information
on  the fibers, as follows:
\begin{stmt}
  \label{stmt:upper-semi-cont}
  For each integer $n$, the set of $x \in X$ such that $\dim_x
  \pi^{-1}(\pi(x)) \ge n$ is closed in $X$.
\end{stmt}
Now suppose that $\A$ is a \emph{local}, Noetherian domain.
Let $\eta \in S$ be the generic point and let $s \in \Spec(\A)$ be the
closed point. 
\begin{stmt}
  \label{stmt:dim-behaviour}
  If $e = \dim \pi^{-1}\eta = \dim X_\eta$ and $f = \dim \pi^{-1} s =
  \dim X_s$, then for each $t \in \Spec(\A)$ we have $f \ge \dim X_t \ge e$.
\end{stmt}

\begin{proof}
  Let $t \in \Spec(\A)$.  Arguing as in \cite{EGA-IV}*{Cor. 13.1.6},
  assertion \ref{stmt:upper-semi-cont} shows that $\dim X_t \ge e$. On
  the other hand, let $\pp \subset \A$ be the prime ideal
  corresponding to $t$.  We may form the fiber product $X
  \times_{\Spec(\A)} \Spec(\A/\pp)$.  The morphism $\Spec(\A/\pp) \to
  \Spec(\A)$ is just the inclusion of the closure of $\{t\}$ in
  $\Spec(\A)$; in particular, $t$ is (identified with) the generic
  point of $\Spec(\A/\pp)$, and $s$ remains the closed point. For any
  point $r \in \Spec(\A/\pp) \subset \Spec(\A)$, the fiber over $r$ of
  $X \times_{\Spec(\A)} \Spec(\A/{\pp})$ identifies with $X_r$. Thus
  the inequality $f \ge \dim X_t$ results from the inequality already
  established.
\end{proof}

\subsection{Existence of sections}
\label{sub:existence-of-sections}
In this section, we consider  a Noetherian, normal domain $\A$.
\begin{stmt}
  \label{stmt:normal-A-section-over-affine}
  Let $Y$ be an affine $\A$-scheme of finite type, and let $x_\eta \in
  Y(k(\eta))$, where $\eta$ is the generic point of $\Spec(\A)$. For
  each prime ideal of height one $\pp \subset \A$, suppose that $x_\pp
  \in Y(\A_\pp)$, and that the image of $x_\pp$ under the natural map
  $Y(\A_\pp) \to Y(k(\eta))$ is $x_\eta$.  Then there is a section $x
  \in Y(\A)$ such that the image of $x$ under the natural map $Y(\A)
  \to Y(\A_\pp)$ is $x_\pp$ for each prime ideal $\pp$ of height one.
\end{stmt}

\begin{proof}
  Write $Y = \Spec(\B)$ for the $\A$-algebra $\B = \A[y_1,\dots,y_n]$.
  The point $x_\eta$ is the same as an $\A$-homomorphism $f_\eta:\B
  \to k(\eta)$. For each prime ideal $\pp \subset \A$ of height one,
  the assumptions mean that $f_\eta(y_i) \in \A_\pp \subset k(\eta)$
  for $1 \le i \le n$, and that the resulting homomorphism $f_\pp:\B \to
  \A_\pp$ determines the point $y_\pp$.

  Since $\A$ is normal, one knows that
  \begin{equation*}
    \A = \cap_{\pp} \A_{\pp}, \quad \text{the intersection taken over all 
      prime ideals $\pp \subset \A$ of height 1};
  \end{equation*}
  see e.g. \cite{liu}*{Lem. 4.1.13}.  We conclude that $f_\eta(y_i)
  \in \A$ for $1 \le i \le n$;  writing $x \in Y(\A)$ for the
  section determined by the resulting homomorphism $f:\B \to \A$, the
  result follows.
\end{proof}

Fix a projective $\A$-scheme $X$ of finite type.  Recall now the
following (see for instance \cite{liu}*{Theorem 3.3.25}):
\begin{stmt}
  \label{stmt:projective-over-dvr-K-points} If $\A$ is a
  discrete valuation ring, then the natural mapping $X(\A) \to
  X(k(\eta))$ is bijective, where $\eta$ is the generic point of
  $\Spec(\A)$.
\end{stmt}

Suppose given an element $x_F \in X(F)$ for each $\A$-algebra $F$ that
is a field. Whenever $\B$ is an $\A$-algebra that is a discrete
valuation ring with field of fractions $F$, write $x_\B \in X(\B)$ for
the section determined by $X(F)$. We make two assumptions:
\begin{enumerate}
\item[(S1)] Whenever the $\A$-algebras $F_1,F_2$ are fields satisfying
  $F_1 \subset F_2$, suppose that $x_{F_2}$ coincides with the image
  of $x_{F_1}$ under the natural map $X(F_1) \to X(F_2)$.
\item[(S2)] Whenever $\B$ is an $\A$-algebra that is a discrete
  valuation ring with field of fractions $F$ and residue field
  $\mathfrak{f}$, we suppose that $x_\mathfrak{f}$ coincides with the
  image of $x_\B$ under the natural map $X(\B) \to X(\mathfrak{f})$.
\end{enumerate}

\begin{stmt}
\label{stmt:section-of-projective-over-normal-A}
  Under the hypotheses (S1) and (S2), there is a unique section $x \in X(\A)$
  such that for each $\A$-algebra $F$ that is a field, the element
  $x_F$ is the image of $x$ under the natural map $X(\A) \to X(F)$.
\end{stmt}

\begin{proof}
  First note that uniqueness of the section $x$ is immediate, e.g.
  since the image of $x$ in $X(k(\eta))$ must coincide with
  $x_{k(\eta)}$.

  We now prove the existence of $x$. In view of the uniqueness of the
  section $x$, it suffices to construct $x$ locally on $\Spec(A)$;
  thus, we may and will suppose that $\A$ is moreover \emph{local}.
  Write $s \in \Spec(\A)$ for the unique closed point, and write $k =
  k(s)$

  Before beginning the proof, choose a very ample invertible sheaf
  $\LL$ on $X$; thus $\LL = i^*\mathcal{O}(1)$ for a suitable closed
  embedding $i:X \to \Proj^n_{/\A}$.  Let $t_0,t_1,\dots,t_n \in
  \LL(X)$ be global sections such that for each $0 \le i \le n$, the
  set $X_{t_i}$ is an affine open $\A$-subscheme, and the affine opens
  $X_{t_i}$ cover $X$. 

  The proof proceeds by induction on $d=\dim \A$. 
  When $d=1$, the domain $\A$ is itself a discrete valuation ring and
  the existence of the desired section $x \in X(\A)$ follows
  immediately from \ref{stmt:projective-over-dvr-K-points}.

  Suppose now that $d>1$ and suppose that the result is true in
  dimension strictly less than $d$.  Let $\pp$ denote a height one
  prime ideal of $\A$. Then the quotient $\A/\pp$ is a Noetherian,
  normal, local domain of dimension $d-1$.  For each $\A/\pp$-algebra
  $F$ which is a field, we have of course the section $x_F$, and it is
  clear that these sections satisfy conditions (S1) and (S2) for the
  ring $\A/\pp$.  Thus, the induction hypothesis gives now a section
  $x' \in X(\A/\pp)$ whose image in $X(F)$ coincides with $x_F$
  for each $\A/\pp$-algebra $F$ that is a field.

  Since $\pp$ has height one and since $\A$ is normal, the
  localization $\A_\pp$ is a discrete valuation ring.  The residue
  field of $\A_\pp$ is $k(\pp)$, the field of fractions of $\A/\pp$.
  By the result when $d=1$, we may find a section $x_\pp \in
  X_{t_j}(\A_\pp)$ whose image in $X(k(\pp))$ coincides with the image
  of $x'$.

  Since $X_{t_j}$ is an affine $\A$-scheme for each $j$, we may now
  apply \ref{stmt:normal-A-section-over-affine} to the restriction to
  $X_{t_j}$ of the sections $\{x_\pp\}$; we then find the required
  section $x \in X(\A)$.
\end{proof}

\section{Reductive groups}
\label{sec:reductive}

We are going to work throughout the remainder of the paper with a
local, normal, Noetherian domain $\A$. Write $K$ for the field of
fractions of $\A$ and $k$ for its residue field. Also, write $\eta \in
\Spec(\A)$ for the generic point; thus $K = k(\eta)$.

\subsection{Group schemes of multiplicative type}
\label{sub:mult}
An $\A$-group scheme $D$ is said to be diagonalizable if there is a
finitely generated Abelian group such that $D \simeq D_\Gamma$, where
$D_\Gamma = \Spec \A[\Gamma]$; here $\A[\Gamma]$ is the group algebra
of $\Gamma$ -- i.e.  the algebra of those $\A$-valued functions on
$\Gamma$ having finite support -- made into a Hopf algebra as usual.
An $\A$-group scheme $M$ is of multiplicative type if it is
diagonalizable locally in the \'etale topology of $\A$; this means
that there is a finite, \'etale extension $\B \supset \A$ such that
$M_{/\B}$ is diagonalizable.

An $\A$-group scheme $T$ is a \emph{torus} if it is a group of
multiplicative type and if locally in the \'etale topology of $\A$ the
group $T$ is of the form $D_\Gamma$ where $\Gamma$ is finitely
generated and free.  The torus $T$ is \emph{split} over $\A$
if it is isomorphic to $D_\Gamma$ as an $\A$-group scheme.

\subsection{Reductive group schemes}
\label{sub:reductive-group-schemes}
Recall that a group scheme $G$ over $\A$ is said to be
\emph{reductive} provided that $G$ is smooth and of finite type over
$\A$ and that the fiber $G_{\bar{t}}$ is a (connected and)
reductive algebraic group for each algebraically closed geometric
point $\bar{t}$ of $\Spec(\A)$. The reductive $G$ is moreover
semisimple if all $G_{\bar{t}}$ are semisimple algebraic
groups.

If $G$ is a group scheme and $T \subset G$ is a subgroup scheme, one
says that $T$ is a \emph{maximal torus} if it is a torus, and if
$T_{t}$ is a maximal torus in $G_t$ for each point $t$ of $\Spec(\A)$.
A result of Grothendieck says:
\begin{stmt}\cite{SGA3}*{Exp. XIV, Cor. 3.20}.
  \label{stmt:reductive-has-max-torus}
  Any reductive group has a maximal torus.
\end{stmt}
The reductive group $G$ is said to be split if it has a split maximal
torus $T$.  If $G$ has a split maximal torus $T$, the root datum of
$G$ with respect to $T$ is  $\mathcal{R}=(X,Y,R,R^\vee)$
where $X = X^*(T)$ is the character group of $T$, $Y=X_*(T)$ is the
group of cocharacters of $T$, $R \subset X$ is the set of roots, and
$R^\vee$ is the set of coroots.  We have the following existence
theorem of Chevalley:
\begin{stmt} \cite{SGA3}*{Exp. XXV, Cor. 1.2} Let $\mathcal{R}$ be
  a root datum. Then there is a split reductive group scheme over $\A$
  with the root datum $\mathcal{R}$.
\end{stmt}

A root datum $\mathcal{R}=(X,Y,R,R^\vee)$ is said to be \emph{of
  adjoint type} if $X = \Z R$. Given any root datum, one constructs
the corresponding adjoint root datum $\mathcal{R}_{\ad}$ and the
morphism $h:\mathcal{R}_{\ad} \to \mathcal{R}$ of root data, as in
\cite{SGA3}*{Exp. XXI, Prop. 6.5.5}. Let $G_{\ad}$ be a split semisimple
$\A$-group scheme with split maximal torus $T'$ and root datum
$\mathcal{R}_{\ad}$.
\begin{stmt}\cite{SGA3}*{Exp. XXIII, Thm. 4.1}
  \label{stmt:adjoint-group}
  There is a unique morphism of $\A$-group schemes $f:G \to G_{\ad}$ which
  defines upon restriction to $T$ a morphism $f_{\mid T}:T \to T'$ and
  induces the mapping $h$ on root data.
\end{stmt}

\subsection{Levi factors}
\label{sub:Levi-factor}

Let $H$ be a smooth and separated \footnote{Note that we will consider
  only group schemes $H$ which are affine over $\A$, and thus
  automatically separated.} group scheme of finite type over $\A$.

Suppose that there exists a closed subgroup scheme $L \subset H$ such
that $L$ is smooth over $\A$ and $L^o$ is reductive. Then we say that
$L$ is a \emph{Levi factor} of $H$ if for each $t \in \Spec(\A)$ the
inclusion $L_t \subset H_t$ induces an isomorphism of
$k(\bar{t})$-group schemes $L_{\bar{t}} \simeq H_{\bar{t}} /
R_u(H_{\bar{t}})$, where $R_u(H_{\bar{t}})$ is the unipotent of
$H_{\bar{t}}$ and $\bar{t}$ is a geometric point over $t$.
Equivalently put: for each $t \in \Spec(\A)$, the subgroup
$L_t$ is a Levi factor of $H_t$ in the usual sense of
linear algebraic groups.

\begin{stmt}
 \label{stmt:levi-extra}
 Let $L \subset H$ be a Levi factor. Write $i:L \to H$ for the
 inclusion, and suppose that there is a homomorphism $\rho:H \to L$ of
 group schemes over $\A$ such that $\rho \circ i =
 \operatorname{id}_{L}$; in other words, $\rho$ is a ``retraction''.
 Let $R = \ker \rho$. If $\rho$ is smooth, then the mapping $\Phi:L
 \ltimes R \to H$ induced by the natural inclusions $L \to H$ and $R
 \to H$ is an isomorphism, where $L \ltimes R$ is the \emph{semidirect
   product} group scheme.
\end{stmt}

\begin{proof}
  The kernel $R = \rho^{-1}(1)$ identifies with the fiber product $H
  \times_{L,\rho} \Spec(\A)$; since $\rho$ is assumed to be smooth,
  $R$ is smooth over $\A$.

  Since $L_t$ is a Levi factor of $H_t$ for each $t$, $\Phi_t$ is an
  isomorphism for each $t \in \Spec(\A)$. We have seen that $R$ is
  smooth, so the group $L \ltimes R$ is smooth. Since both $L \ltimes
  R$ and $H$ are smooth -- hence flat -- over $\A$, it now follows
  from \cite{SGA1}*{Exp. I, Prop.  5.7} that $\Phi$ is itself an
  isomorphism.
\end{proof}

\subsection{The identity component}

Let again $H$ be a smooth and separated group scheme of finite type
over $\A$.
\begin{stmt}
  There is a smooth, normal, and open subgroup scheme $H^o \subset H$
  which is the union of the connected components of the groups $H_t$
  for $t \in \Spec(\A)$.
\end{stmt}

For us, an important property of the connected component is the
following:
\begin{stmt}
  \label{stmt:reductive-o-closed}
  If $H^o$ is reductive, then $H^o$ is closed in $H$.
\end{stmt}

\begin{proof}
  Since $H^o$ is reductive, there is a maximal torus $T \subset H^o$.
  Now, it follows from \cite{SGA3}*{Exp. XIX, Thm. 2.5} that the
  Weyl group $W= N_{H^o}(T)/C_{H^o}(T) = N_{H^o}(T)/T$ is (represented
  by) a finite \'etale group scheme over $\A$. Since $H$ is assumed to
  be a separated group scheme over $\A$, it follows from
  \cite{SGA3}*{Exp. XVI, Cor. 1.4} that the inclusion $H^o \subset H$ is
  a closed immersion, as required.
\end{proof}

\subsection{Central isogenies}
Let $G$ and $G'$ be reductive groups over $\A$. An $\A$-homomorphism
$f:G \to G'$ is said to be a central isogeny if $f$ is faithfully
flat, finite and if $\ker f$ is a central subgroup of $G$.  If $f$ is
an \'etale central isogeny, observe that $f_s$ and $f_\eta$ are
separable central isogenies in the usual sense of an algebraic group
over a field.

\begin{stmt}
  Let $f:G \to G'$ be an \'etale central isogeny. Then $\ker f$ is a
  closed and central subgroup scheme of $G$ which is finite, \'etale
  and of multiplicative type over $\A$.
\end{stmt}

\begin{proof}
  Let $D = \ker f$. Since $f:G \to G'$ is finite and \'etale, upon
  base change we see that $D=\ker f \to \Spec(\A)$ is finite and
  \'etale as well. To see that $D$ is of multiplicative type, we may
  replace $\A$ by a finite, \'etale, local extension and show that
  $D$ is diagonalizable. Thus, using \cite{SGA3}*{Exp. XXII, Cor.
    4.2.13} we may suppose that there are split maximal tori $T
  \subset G$ and $T' \subset G'$ such that $f_{\mid T}$ factors as a
  morphism $f_{\mid T}:T \to T'$.  Since $\ker f$ is \'etale over
  $\A$, it follows from \cite{SGA3}*{Exp. XXII, Cor. 4.2.8} that
  $\ker f = \ker f_{\mid T}$, so that $\ker f$ is a closed subgroup of
  the torus $T$. It now follows from \cite{SGA3}*{Exp. VII, Cor.
    3.4} that $\ker f$ is diagonalizable, as required.
\end{proof}



\subsection{Some centralizers}
\label{sub:centralizer}

Let $H$ be a group scheme which is smooth of finite type over $\A$.
Let $D \subset H$ be a subgroup scheme of multiplicative type.

\begin{stmt}{\cite{SGA3}*{Exp. XI, Cor. 5.3}}
  \label{stmt:mult-centralizer}
  The centralizer $C_H(D)$ and the normalizer $N_H(D)$ are closed subgroup
  schemes which are smooth over $\A$.
\end{stmt}

Now suppose that $G$ is a reductive $\A$-group scheme and that $D
\subset H$ is a \emph{smooth} subgroup scheme of multiplicative type.
\begin{stmt}
  \label{stmt:mult-centralizer-in-G}
  The centralizer $C_G(D)$ is a closed and smooth subgroup scheme
  whose identity component $C_G(D)^o$ is a reductive $\A$-group
  scheme.
\end{stmt}

\begin{proof}[Sketch]
  Since $C_G(D)$ is closed and smooth over $\A$ by
  \ref{stmt:mult-centralizer}, it suffices to show that
  $C_G(D)_{\bar{t}}$ has reductive identity component for each $t \in
  \Spec(\A)$; thus, it is enough to prove the result when $G$ is
  reductive over an algebraically closed field.  In that case, $D$ is
  diagonalizable, and arguing by induction on $\dim G$ one quickly
  reduces to two cases: $D$ a torus, in which case $C_G(D)$ is a Levi
  factor of a parabolic subgroup of $G$, and $D$ a cyclic group of
  order invertible in the field, in which case the result follows from
  \cite{steinberg-endomorphisms}*{Cor. 9.3}.
\end{proof}

\begin{rem}
  In fact, the preceding result remains valid for \emph{any}
  diagonalizable group $D$. It seems to be difficult to find a
  reference for this more general fact. In case $D$ is smooth, the
  assertion that $C_G(D)^o$ is reductive may also be deduced from a
  result of Richardson \cite{richardson}*{Prop.  10.1.5}; I thank
  Gerhard R\"ohrle for pointing out this reference to me.
\end{rem}

\begin{stmt}
  \label{stmt:D-contained-in-T}
  There is a maximal torus of $G$ centralized by $D$. Moreover,
  if $T$ is any maximal torus of $G$ which is centralized by $D$, then $D
  \subset T$.
\end{stmt}

\begin{proof}
  Let $M = C_G(D)$ be the centralizer of $D$ in $M$. Then by
  \ref{stmt:mult-centralizer-in-G} the identity component $M^o$ is a
  reductive subgroup scheme, hence $M^o$ contains a maximal torus $T$
  \ref{stmt:reductive-has-max-torus}. It follows from
  \cite{SGA3}*{Exp. XII, Prop. 1.17} that $T$ is a maximal torus in
  $G$ as well.  

  Since $D$ is central in $M$, \cite{SGA3}*{Exp. XII, Lem. 4.5}
  shows for any maximal torus $T$ of $M^o$ that the inclusion $D
  \subset M$ factors through $T$; i.e. $D$ is a subgroup scheme of $T$
  as required.
\end{proof}

\begin{stmt}
  Let $T \subset G$ be a torus. Then $C_G(T)$ is a closed and
  reductive subgroup scheme.
\end{stmt}

\begin{proof}
  In view of \ref{stmt:mult-centralizer-in-G}, it only remains to show
  that $C_G(T) = C_G(T)^o$, i.e. that $C_G(T)_x$ is connected \ for
  every $x \in \Spec(\A)$; that  requirement follows e.g. from
  \cite{SGA3}*{Exp. XIX, \S 1.3} or from
  \cite{springer-LAG}*{Theorem 6.4.7}.
\end{proof}

\subsection{The derived group}
\label{sub:derived-group}
Let $\Der(G)$ be the derived subgroup scheme of $G$. Then:
\begin{stmt}
  \cite{SGA3}*{Exp. XXII, Thm. 6.2.1}
  $\Der(G)$ is a closed, normal subgroup scheme which is smooth and
  semisimple over $\A$. 
\end{stmt}

\begin{stmt}
  \label{stmt:max-torus-of-Der(G)}
  \cite{SGA3}*{Exp. XXII, Prop. 6.2.7}
  If $T$ is a split maximal torus of $G$, there is a maximal torus $T' 
  \subset \Der(G)$ contained in $T$.
\end{stmt}


\subsection{Good and very good primes}
\label{sub:good-primes}

Let $\FFF$ denote an arbitrary field, and
let $H$ be a geometrically quasisimple algebraic group over $\FFF$ with
absolute root system \footnote{The absolute root system of $G$ is the
  root system of $G_{/\FFF_\sep}$ where $\FFF_\sep$ is a separable closure of
  $\FFF$.}  $R$. The characteristic $p$ of $\FFF$ is said to be a bad prime
for $R$ in the following circumstances: $p=2$ is bad whenever $R \not
= A_r$, $p=3$ is bad if $R = G_2,F_4,E_r$, and $p=5$ is bad if
$R=E_8$.  Otherwise, $p$ is good.  [Here is a more intrinsic
definition of good prime: $p$ is good just in case it divides no
coeficient of the highest root in $R$].

If $p$ is good, then $p$ is said to be very good provided that either
$R$ is not of type $A_r$, or that $R=A_r$ and $r \not
\congruent -1 \pmod p$.

If $H$ is reductive, one may apply \cite{KMRT}*{Theorems 26.7 and
  26.8} \footnote{\cite{KMRT} only deals with the semisimple case; the
  extension to a general reductive group is not difficult to handle,
  and an argument is sketched in the footnote found in
  \cite{mcninch-testerman}*{\S2.4}.}  to see that there is a possibly
inseparable isogeny
\begin{equation}
  \label{eq:cover-G}
   \quad T \times \prod_{i=1}^r H_i \to H
\end{equation}
for some $\FFF$-torus $T$ and some $r \ge 1$, where for $1 \le i \le r$
there is an isomorphism $H_i \simeq R_{E_i/\FFF} J_i$ for a finite
separable field extension $E_i/\FFF$ and a geometrically quasisimple,
simply connected $E_i$-group scheme $J_i$; here, $R_{E_i/\FFF} J_i$
denotes the ``Weil restriction'' of $J_i$ to $\FFF$.

Then $p$ is good, respectively very good, for $H$ if and only if that
is so for $J_i$ for every $1 \le i \le r$. Let $\FFF_\alg$ be an algebraic
closure of $F$. Since the groups $J_{i/\FFF_\alg}$ are uniquely determined
by $H_{/\FFF_\alg}$ up to central isogeny, the notions of good and very
good primes depend only on the group $H_{/\FFF_\alg}$, and these notions
depend only on the central isogeny class of the derived group of
$H_{/\FFF_\alg}$

Let now $G$ be a \emph{split semisimple} scheme over $\A$ with split
maximal torus $T$ and corresponding root datum $\mathcal{R}$. Suppose
that the characteristic of the residue field $k$ is very good for $G_s
= G_{/k}$. For any point $t$ of $\Spec(\A)$, either the characteristic
of $k(t)$ is zero, or is the same as the characteristic of $k$; thus
also the characteristic of $k(t)$ is very good for $G_t$.

\begin{stmt}
  \label{stmt:very-good-p-adjoint-group}
  Let $f:G \to G_{\ad}$ be the morphism of \ref{stmt:adjoint-group},
  where $G_{\ad}$ is the $\A$-group scheme with the root datum
  $\mathcal{R}_{\ad}$. Then $f$ is an \'etale central isogeny.
\end{stmt}
\begin{proof}
  That $f$ is a central isogeny follows from \cite{SGA3}*{Exp. XXII,
    Prop. 4.2.10}. To see that $f$ is \'etale, it suffices by
  \cite{SGA1}*{Exp. I, Cor. 5.9} to observe that the mapping
  $f_t:G_{t} \to G_{\ad,t}$ is \'etale for each $t$ in $\Spec(\A)$.
  In view of our assumptions on the characteristic, one may use the
  descriptions found in \cite{Hum95}*{0.13} to see that the tangent
  mapping $df_{t}$ must be an isomorphism, whence the
  required assertions.
\end{proof}

\subsection{Strongly standard reductive group.}
\label{sub:strongly-standard} 

Consider reductive groups $H$ over $\A$ 
which are direct products
\begin{equation*}
  (*) \quad
  H=H_1  \times_S T
\end{equation*}
where $T$ is a torus over $\A$, and where $H_1$ is a semisimple group
over $\A$ such that the characteristic of $K$ is very good for
$H_{1/K}$ and such that the characteristic of $k$ is very good for
$H_{1/k}$.

Let $G$ be a reductive $\A$-group. Then $G$ will be said to be
\emph{$D$-standard} if there exists a reductive group $H$ of the form
$(*)$, an $\A$-subgroup $D \subset H$ of multiplicative type, and an
\'etale $\A$-isogeny between $G$ and the $\A$-group $C_H(D)^o$
\footnote{If there is such an isogeny, then $C_H(D)^o$ is of course
  reductive, so the definition is independent of
  \ref{stmt:mult-centralizer-in-G}; in particular, there is no need to
  insist in the definition that $D$ be smooth.}.

Similarly, $G$ will be said to be \emph{$T$-standard} if there
exists a reductive group $H$ of the form $(*)$, an $\A$-torus $T
\subset H$, and an \'etale $\A$-isogeny between
$G$ and the smooth reductive $\A$-group $C_H(T)$.

Of course, any $T$-standard group is $D$-standard.

\begin{rem}
  Let $\FFF$ be a field.  In \cite{mcninch-optimal}, the term
  \emph{strongly standard} reductive group was used for what we call here a
  $T$-standard group scheme over $\FFF$. In \cite{mcninch-testerman}, the
  term \emph{strongly standard} reductive group was used for what we call here a
  $D$-standard group scheme over $\FFF$.
\end{rem}

We have evidently
\begin{stmt}
  If $G$ is $D$-standard, respectively $T$-standard, then for $t \in
  \Spec(\A)$, the fiber $G_t$ is $D$-standard, respectively $T$-standard.
\end{stmt}

\begin{stmt}
  \label{stmt:smooth-centralizers}
  Let $G$ be a $D$-standard reductive group over a field $\FFF$,
  let $g \in G(\FFF)$ and $X \in \glie(\FFF)$. Then the centralizers
  $C_G(g)$ and $C_G(X)$ are smooth $\FFF$-subgroup schemes of $G$.
\end{stmt}

\begin{proof}
  \cite{mcninch-testerman}*{Prop. 12}.
\end{proof}

\begin{stmt}
  Let $\LL$ be a free $\A$-module of finite rank $n$. Then the
  reductive $\A$-group $\GL(\LL)$ is $T$-standard (hence also
  $D$-standard).
\end{stmt}
\begin{proof}[Sketch]
  Indeed, let $p$ denote the
  characteristic of the residue field $k$. If $p=0$, nothing needs to
  be said, so assume $p>0$. If $n \not \equiv 0 \pmod p$ then
  $H=\SL(\LL) \times \Gm$ is of the form $(*)$, and multiplication
  defines an \'etale isogeny $H \to \GL(\LL)$.
  
  If $n \equiv 0 \pmod p$, then $H = \SL(\LL \oplus \A)$ has the form
  $(*)$, and $\GL(\LL)$ is isomorphic to the centralizer in $H$ of a
  suitable split $\A$-torus $S \subset H$.
\end{proof}

\subsection{Parabolic subgroups}
\label{sub:parabolics}

Let $G$ be a reductive group scheme over $\A$, and let $P \subset G$
be an $\A$-subgroup scheme. One says that $P$ is an
\emph{$\A$-parabolic subgroup scheme} of $G$ if $P$ is smooth over
$\A$ and if $P_t$ is a parabolic subgroup of $G_t$ for each point $t$
of $\Spec(\A)$.

We recall the following:
\begin{stmt}
  \label{stmt:parabolic-scheme-projective}
  \cite{SGA3}*{Exp. XXVI, Cor. 3.5}.
  Consider the functor $\Par$ defined
  for  commutative $\A$-algebras $\Lambda$ by the rule
  \begin{equation*}
    \Par(\Lambda)  = \text{set of all $\Lambda$-parabolic subgroup schemes of $G_{/\Lambda}$}.
  \end{equation*}
  Then $\Par$ is (represented by) a smooth and projective scheme over $\Spec(\A)$.
\end{stmt}

\begin{stmt}
  \label{stmt:parabolic-determined-by-P_eta}
  Let $P,Q \subset G$ be $\A$-parabolic subgroup schemes, and write
  $\eta$ for the generic point of $\Spec(\A)$.
  \begin{enumerate}
  \item[(a)] If $P_\eta = Q_\eta$, then $P = Q$.
  \item[(b)] If $P_\eta$ and $Q_\eta$ are conjugate by an element of
    $G(k(\eta))$, then $P = \Int(g)Q$ for a section $g \in
    G(\A)$.
  \end{enumerate}
\end{stmt}

\begin{proof}
  For (a), note first that by \ref{stmt:parabolic-scheme-projective}
  the scheme $\Par$ of parabolic subgroups of $G$ is projective --
  hence in particular, separated -- over $\A$.  Thus if the
  restrictions of two $\A$-sections of $\Par$ to the dense open subset
  $\{\eta\} \subset \Spec(\A)$ coincide, then the sections coincide by
  \cite{liu}*{Prop.  3.3.11}.

  For (b), let $\beta(\operatorname{Dyn}(G))$ be the scheme of
  \emph{types} of the parabolic subgroup schemes of $G$, and for an
  $\A$-parabolic subgroup scheme $P$ of $G$, we will write
  $\underline{t}(P) \in \beta(\operatorname{Dyn}(G))(\A)$ for the type
  of $P$; cf.  \cite{SGA3}*{Exp. XXVI, Defn 3.4}.  Thus
  $\underline{t}:\Par \to \beta(\operatorname{Dyn}(G))$ is a morphism
  of schemes. Since $\{\eta\}$ is dense in $\Spec(\A)$ and since
  $P_\eta$ and $Q_\eta$ are $G(k(\eta))$-conjugate, it follows that
  $\underline{t}(Q) = \underline{t}(P)$. Since $\A$ is local,
  assertion (b) now follows from \cite{SGA3}*{Exp. XXVI, Cor. 5.5}.
\end{proof}

\begin{stmt}
  \label{stmt:parabolic-over-dvr}
  Assume that $\A$ has dimension 1; i.e. assume that $\A$ is a discrete
  valuation ring.  If $Q \subset G_\eta= G_{/K}$ is a $K$-parabolic
  subgroup, there is a unique parabolic $\A$-subgroup scheme $P
  \subset G$ inducing $Q$ on base-change -- i.e. $Q = P_\eta$.
\end{stmt}

\begin{proof}
  Indeed, uniqueness follows from
  \ref{stmt:parabolic-determined-by-P_eta}(a). Since the scheme $\Par$
  of parabolic subgroups of $G$ is projective
  \ref{stmt:parabolic-scheme-projective}, and since $\A$ is a discrete
  valuation ring, it follows that its $K$-points are the same as its
  $\A$-points \cite{liu}*{Theorem 3.3.25}.
\end{proof}

The validity of \ref{stmt:parabolic-over-dvr} indeed requires some
hypothesis on $\A$. Notice that if $G = \GL_2$, then the scheme of
Borel subgroups of $G$ identifies with the projective line
$\Proj^1_\A$.  If $k$ is a field and $\A$ is the 2 dimensional
(regular, hence normal) local domain $\A=k[x,y]_{(x,y)}$, then the
$K$-point $(x:y) \in \Proj^1_\A(K)$ does not determine an $\A$-section
of $\Proj^1_\A$.

\subsection{Cocharacters and parabolic subgroups}
\label{sub:cochar-and-parabolic-over-A}

If $H$ is an algebraic group over a field $\FFF$, a cocharacter of $H$
is an $\FFF$-homomorphism $\Gm \to H$. In this paper, we will be
interested more generally in homomorphisms of group schemes from the
multiplicative group to a given group. 

Let $H$ be an $\A$-group scheme, consider a representation of $H$ on a
free $\A$-module of finite rank $V$ given by the comodule map $\rho:V
\to V \tensor_\A \A[H]$.  If $\phi:\Gm \to H$ is an $\A$-homomorphism
of group schemes, one obtains a representation of $\Gm$ on $V$ with
co-module map $(1 \tensor \phi^*) \circ \rho:V \to V \tensor_\A
\A[\Gm] = V \tensor_\A \A[t,t^{-1}]$.  One now finds a direct sum
decomposition $V = \bigoplus_{n \in \Z} V(\phi;n)$ where
$V(\phi;n)$ is the $n$-weight space; i.e. 
\begin{equation*}
  V(\phi;n) = \{v \in V \mid 
  (1 \tensor \phi^*) \circ \rho(v) = v \tensor t^n\}.
\end{equation*}
We apply this especially when $H$ is smooth and of finite type over
$\A$, so that $\hlie = \Lie(H)$ is a free $\A$-module of finite rank
on which $H$ acts by the adjoint representation.

As to the \emph{existence} of $\A$-homomorphisms $\Gm \to H$,
we first note that since the domain $\A$ is assumed to be
\emph{normal}, we have:
\begin{stmt}[\cite{SGA3}*{Exp. X, Lem. 8.4}]
  Let $D,H$ be group schemes over $\A$, where $D$ is of multiplicative
  type, and $H$ is smooth. Let $\phi:H_\eta \to D_\eta$ be a
  homomorphism of group schemes over $K = k(\eta)$. Then there is a
  unique homomorphism $\psi:H \to D$ of group schemes over $\A$ such
  that $\phi = \psi_\eta$.
\end{stmt}

An immediate consequence is the following:
\begin{stmt}
  \label{stmt:cochar-over-A}
  Let $T$ be an $\A$-torus, and let $\phi:\Gm_{,\eta} \to T_{\eta}$ be a
  cocharacter over $k(\eta)=K$.  Then there is a unique $\A$-homomorphism of
  group schemes $\psi:\Gm \to T$ such that $\phi = \psi_{\eta}$.
\end{stmt}

Let again $G$ be a reductive group scheme over $\A$.  Suppose now that
$\phi:\Gm \to G$ is an $\A$-homomorphism of group schemes. Composing
$\phi$ with the (left) regular representation of $G$, the algebra
$\A[G]$ becomes a locally finite representation of $\Gm$; let us write
it as the direct sum of its weight spaces
\begin{equation*}
  \bigoplus_{n \in \Z} \A[G]_n.
\end{equation*}
Form the ideal $I$ generated by $\sum_{n > 0} \A[G]_n$.  Then $P(\phi) = P_G(\phi)
= \Spec(\A[G]/I)$ is a closed subgroup scheme of $G$.

\begin{stmt}
  \label{stmt:parabolic-P(phi)}
  $P(\phi) = P_G(\phi)$ is a parabolic subgroup scheme of $G$
  with $\Lie P(\phi) = \bigoplus_{n \ge 0} \glie(\phi;n).$
\end{stmt}

\begin{proof}
  It follows from \cite{springer-LAG}*{Prop. 8.4.5 and Theorem 13.4.2}
  that $P=P(\phi)$ determines a parabolic subgroup of each fiber upon
  base-change; in particular, $P_t$ is smooth over $k(t)$ for each $t
  \in \Spec(\A)$. To see that $P$ is smooth over $\A$, we may first
  replace $\A$ by an \'etale local extension and thus suppose the
  image of $\phi$ to lie in a split maximal torus of $G$. Then $P$ is
  a standard parabolic and hence smooth.
\end{proof}

Let $P$ be any parabolic subgroup scheme of $G$. 
\begin{stmt}{\cite{SGA3}*{Exp. XXVI, Prop. 1.6}}
  \label{stmt:Ru(P)}
  There is a largest normal subgroup scheme $R=R_u(P) \subset P$ which
  is smooth over $\A$ and has connected and unipotent geometric
  fibers. The geometric fiber $R_{\bar{t}}$ is the unipotent radical
  of $P_{\bar{t}}$ for each  $t \in \Spec(\A)$. If $P =
  P(\phi)$ for an $\A$-homomorphism $\phi:\Gm \to G$, then $\Lie
  R_u(P) = \bigoplus_{n > 0} \glie(\phi;n)$.
\end{stmt}

Recall from \S \ref{sub:Levi-factor} that by a Levi subgroup scheme of
$P$ we mean a closed and smooth subgroup scheme $L \subset P$ such
that $L^o$ is reductive and such that $L_t$ is a Levi factor of $P_t$
for each $t \in \Spec(\A)$.  Since
$\A$ is assumed to be local, we have:
\begin{stmt}{\cite{SGA3}*{Exp. XXVI, Cor. 2.3, 2.4}}
  \label{stmt:parabolic-has-levi}
  $P$ contains a Levi subgroup scheme, and $P$ contains a maximal torus.
\end{stmt}

Using \cite{SGA3}*{Exp. XXVI, Prop. 1.6}, we see that the
conditions of \ref{stmt:levi-extra} hold; thus $P$ is isomorphic to
the semidirect product $L \ltimes R_u(P)$ for any Levi factor $L$ of
$P$.

In fact, we can be a bit more precise regarding Levi subgroups and
maximal tori, as follows:
\begin{stmt}
  \label{stmt:Levi-of-P-via-Phi}
  If $P = P(\phi)$ for some $\A$-homomorphism $\phi:\Gm \to G$, then
  the centralizer in $G$ of the image of $\phi$ is a Levi subgroup
  scheme of $P$.
\end{stmt}

\begin{proof}
  Write $L$ for the centralizer of the image of $\Phi$.  Then $L$ is a
  closed subgroup scheme of $P$, and according to
  \ref{stmt:mult-centralizer-in-G}, $L$ is smooth over $\A$. Using
  \cite{springer-LAG}*{Theorem 13.4.2} we see that $L$ is indeed a
  Levi subgroup scheme of $P$.
\end{proof}

\begin{stmt}
  \label{stmt:P-is-P(phi)}
  If $T \subset P$ is a maximal torus, then $P = P(\phi)$ for some
  $\A$-homomorphism $\phi:\Gm \to T$. In particular, there is a Levi
  subgroup scheme $L \subset P$ which contains $T$.
\end{stmt}

\begin{proof}
  If $\eta$ denotes the generic point of $\Spec(\A)$, one knows that
  $P_\eta$ is the parabolic subgroup determined by some cocharacter
  $\phi_0$ of the maximal torus $T_\eta \subset P_\eta$ (see e.g.
  \cite{mcninch-optimal}*{Lem. 6}). Since $T$ is an $\A$-torus and
  since $\A$ is normal, use \ref{stmt:cochar-over-A} to find an
  $\A$-homomorphism $\phi:\Gm \to T$ such that $\phi_0 = \phi_{\eta}$.
  It follows from \ref{stmt:parabolic-determined-by-P_eta}(a) that $P =
  P(\phi)$. Finally, \ref{stmt:Levi-of-P-via-Phi} gives the required
  Levi subgroup scheme of $P$.
\end{proof}

\section{Nilpotent elements and the instability parabolic over a field}
\label{sec:instab-parabolic-field}

In this section, we let $\FFF$ be an arbitrary field, and we suppose that $G$
is a $D$-standard reductive group over $\FFF$ with Lie algebra $\glie$.
Let $X \in \glie(\FFF)$ be a nilpotent element.

\subsection{Associated cocharacters}
\label{sub:associated}

We write $X_*(G)$ for the collection of $\FFF$-homomorphisms $\Gm \to G$.
If $\Psi \in X_*(G)$, recall that -- as in
\ref{sub:cochar-and-parabolic-over-A} -- we may write
 $\glie = \bigoplus_{n \in \Z}
\glie(\Psi;n)$ where we regard the Lie algebra $\glie$ as a $G$-module
via the adjoint representation.

A cocharacter $\Psi \in X_*(G)$ is said to be associated with $X$ (see
\cite{jantzen-nil}*{\S 5}) if the following conditions hold:

\begin{enumerate}
\item[(A1)] $X \in \glie(\Psi;2)$, and
\item[(A2)] there is a maximal torus $S$ of $C_G(X)$ such that
  $\Psi \in X_*(L_1)$ where $L=C_G(S)$ and $L_1 = (L,L)$ is its derived group.
\end{enumerate}

By regarding the nilpotent element $X$ as an \emph{unstable} vector in
the $G$-representation $\glie$ and using the notions of optimal
cocharacters and the instability parabolic due to Kempf and to
Rousseau, one finds:
\begin{stmt}
  \label{stmt:assoc-cochar}
  Let $X \in \glie$ be nilpotent. 
  \begin{enumerate}
  \item[(a)] There is a cocharacter $\Psi$ associated with $X$.
  \item[(b)] If $\Psi$ is associated to $X$, then $C_G(X) \subset P(\Psi)$.
  \item[(c)] The unipotent radical of $C=C_G(X)$ is defined over $F$,
    and is an $\FFF$-\emph{split} unipotent group.  
  \item[(d)] If the cocharacter $\Psi$ is associated with $X$, then
    $L=C \cap C_G(\Psi(\G_m))$ is a Levi factor of $C$.
  \item[(e)] If $\Psi,\Phi \in X_*(G)$  are associated
    with $X$, then $\Psi = \Int(x) \circ \Phi$ for a unique $x \in U(K)$.
  \item[(f)] The parabolic subgroups $P(\Psi)$ for cocharacters $\Psi$
    associated with $X$ all coincide.
  \end{enumerate}
\end{stmt}

\begin{proof}
  In the ``geometric case'' -- when $\FFF$ is algebraically closed -- (a)
  is essentially a consequence of Pommerening's -- and more recently,
  Premet's -- proof of the Bala-Carter theorem; Premet's proof
  \cite{premet} avoids case analysis and uses results of geometric
  invariant theory due to Kempf and Rousseau.  Working over any ground
  field $\FFF$, (a) is in \cite{mcninch-rat}*{Theorem 26}. Now (b) is
  \cite{jantzen-nil}*{Prop.  5.9}. (c) and (e) follow from parts (3)
  and (4) of \cite{mcninch-optimal}*{Prop/Defn 21}. (d) is essentially
  a consequence of results in \cite{premet}; see \cite{mcninch-rat}*{Cor. 20 and
    Cor. 29}. Finally, (f) is \cite{mcninch-optimal}*{Prop/Defn
    21(5)}.
\end{proof}

We write $P_X$ for the common parabolic subgroup of part (f) of
\ref{stmt:assoc-cochar}; we say that $P_X$ is the \emph{instability
  parabolic subgroup} attached to $X$.

\subsection{The stabilizer of the line through $X$}
Let $N_G(X) \subset G$ be the stabilizer of the line $[X] \in
\Proj(\glie)$, where $\Proj(\glie)$ denotes the projective $\FFF$-variety
formed from the vector space $\glie$.  Then:
\begin{stmt}{\cite{mcninch-rat}*{Lem. 23}}
   $N_G(X)$ is a smooth $\FFF$-subgroup of $G$.
\end{stmt}
Of course, any cocharacter of $G$ associated to $X$ is a cocharacter
of $N_G(X)$. A more precise version of \ref{stmt:assoc-cochar}(a) is
as follows:
\begin{stmt}{\cite{mcninch-rat}*{Lem. 25}}
  Let $S$ be a maximal torus of $N_G(X)$. Then there is a unique
  cocharacter of $S$ which is associated to $X$.
\end{stmt}

\subsection{Almost associated cocharacters}

Let $P_X$ be the instability parabolic subgroup attached to $X$, let
$\Phi$ be a cocharacter of $G$, and let $\FFF_\sep$ be a separable closure
of $\FFF$. We say that $\Phi$ is \emph{almost associated} to $X$ provided
that $\Int(g) \circ \Phi$ is a cocharacter of $G_{/\FFF_\sep}$ associated
to $X$ for some $g \in P_X(\FFF_\sep)$.

\begin{stmt}
  \label{stmt:almost-assoc-exists}
  Let $S \subset P_X$ be a maximal torus. Then there is a unique
  cocharacter $\Phi$ of $S$ which is almost associated to $X$.  The
  cocharacter $\Phi$ is associated to $X$ if and only if $S$ contains
  a maximal torus of $N_G(X)$.
\end{stmt}

\begin{proof}
  For the existence of $\Phi$, let $S_1$ be a maximal torus of $N_G(X)$, and let
  $S_0$ be a maximal torus of $P$ containing $S_1$. Then $S$ and $S_1$
  are maximal tori of $P$ and hence are conjugate by an element $g\in
  P(\FFF_\sep)$. If $\Psi$ is the cocharacter of $S_1$ associated to $X$,
  then $\Int(g) \circ \Psi$ is a cocharacter of $S$ which is almost
  associated to $X$, as required.

  We now argue the uniqueness. Since $\Phi$ is $P(\FFF_\sep)$-conjugate to
  a cocharacter associated to $X$, one knows that $\Phi$ is an
  \emph{optimal} cocharacter for the unstable vector $X$ in the sense
  of geometric invariant theory; cf.  \cite{mcninch-rat}*{\S3}. Thus
  the unicity is a consequence of the result of Kempf and of Rousseau;
  cf. \cite{mcninch-rat}*{Prop. 13(4)}.

  The remaining assertion is clear.
\end{proof}

\begin{stmt}
  \label{stmt:almost-assoc-subspace}
  If the cocharacters $\Phi,\Psi$ are almost associated to $X$, then
  \begin{equation*}
    \sum_{j \ge 2} \glie(\Psi;j) = \sum_{j \ge 2} \glie(\Phi;j).
  \end{equation*}
\end{stmt}

\begin{proof}
  Indeed, we have $\Phi = \Int(g) \circ \Psi$ for some $g \in
  P(\FFF_\sep)$, so the assertion follows from the fact that $\sum_{j \ge
    2} \glie(\Psi;j)$ is $\Ad(P)$-stable.
\end{proof}

We have:
\begin{stmt}
  \label{stmt:dense-P-orbit}
  Let $\Psi$ be a cocharacter of $G$ which is almost associated to
  $X$. 
  \begin{enumerate}
  \item[(a)] The $\Ad(P_X)$-orbit of $X$ is dense in $\sum_{j \ge 2}
    \glie(\Psi;j)$.
  \item[(b)] Write $X = \sum_{j \ge 2} X_j$ with $X_j\in
    \glie(\Psi;j)$.  Then $X$ is $\Ad(P)(\FFF)$-conjugate to $X_2$, and
    $\Psi$ is a cocharacter associated to $X_2$.
  \end{enumerate}
\end{stmt}

\begin{proof}
  (a) follows by combining \ref{stmt:almost-assoc-subspace} with 
  \cite{jantzen-nil}*{Prop. 5.9(c)}.
  
  The conjugacy statement of (b) follows from
  \cite{mcninch-rat}*{Prop. 34}. It is then clear that $\Psi$ is
  almost associated to $X_2$. Since $\Psi$ is a cocharacter of a
  maximal torus of $N_{X_2}$, it follows from
  \ref{stmt:almost-assoc-exists} that $\Phi$ is associated to $X_2$.
\end{proof}

In particular, \ref{stmt:dense-P-orbit}(b) implies:
\begin{stmt}
  If the cocharacter $\Psi$ is almost associated with $X$ and if $X
  \in \glie(\Psi;2)$, then $\Psi$ is associated with $X$.
\end{stmt}

\subsection{The Bala-Carter theorem}
\label{sub:bala-carter}

For a $D$-standard reductive group over a field $\FFF$, the geometric
nilpotent orbits -- i.e. the nilpotent orbits of $G_{/\FFF_\sep}$ -- are
described by the Bala-Carter theorem. Let us suppose that $\FFF = \FFF_\sep$.

Recall that a parabolic subgroup $P \subset G$ is distinguished if 
\begin{equation*}
  \dim P/U = \dim U/(U,U) + \dim Z
\end{equation*}
where $U$ is the unipotent radical of $P$, and $Z$ is the center of
$G$. A nilpotent element $X \in \glie$ is said to be distinguished if
a maximal torus of $C=C_G(X)$ is central in $C$; if $X$ is
distinguished, then the instability parabolic subgroup $P_X$ is
distinguished.

Each parabolic subgroup has an open orbit -- known as the Richardson
orbit -- on $\Lie R_uP$; any element of this orbit is known as a
Richardson element for $P$.

We have the following important theorem:
\begin{stmt}(The Bala-Carter theorem) 
  \label{stmt:Bala-Carter-theorem}
  Let $L$ be a Levi factor of a parabolic subgroup of $G$, and let $P
  \subset L$ be a distinguished parabolic subgroup of $L$. The map
  which associates to $(L,P)$ the $G$-orbit of a Richardson element
  for $P$ determines a bijection between the $G$-orbits of such pairs
  $(L,P)$ and the $G$-orbits on nilpotent elements in $\Lie(G)$.
\end{stmt}

This theorem was originally proved by Bala and Carter in the case
where $p$ is ``very large''. Pommerening gave a proof in good
characteristic, using some case analysis in a few situations. Premet
\cite{premet} gave recently a short and conceptual proof of this
theorem. See also \cite{jantzen-nil}*{\S4}.

If the orbit of a nilpotent element $X \in \Lie(G)$ corresponds via
the Bala-Carter theorem to the pair $(L,P)$, then the $G$-orbit of
$(L,P)$ -- or, abusing terminology somewhat, just the pair $(L,P)$ --
is said to be the Bala-Carter datum for $X$.

\begin{stmt}
  \label{stmt:choose-favorite-fields}
  Let $\FFF$ and $\FFF'$ be algebraically closed fields of the same
  characteristic, suppose that $G$ and $G'$ are $D$-standard reductive
  groups respectively over $\FFF$ and $\FFF'$ with identical root data, let
  $X \in \Lie(G)$ and $X' \in \Lie(G')$ be nilpotent elements with the
  same Bala-Carter datum, and let $C,C'$ be their respective
  centralizers. Then:
  \begin{enumerate}
  \item the root datum of a Levi factor of $C$ identifies with the
    root datum of a Levi factor of $C'$, and
  \item $C/C^o \simeq C'/C'^o$.
  \end{enumerate}
\end{stmt}

\begin{proof}
  Indeed, we may choose an algebraically closed field $\FFF''$ containing
  both $\FFF$ and $\FFF'$. We thus see that it suffices to prove the
  result when $\FFF \subset \FFF'$.
  
  But then the Bala-Carter theorem implies that $X$ and $X'$ are
  conjugate by an element of $G'$, and the result is immediate.
\end{proof}

\section{Nilpotent sections and the instability parabolic over $\A$}
\label{sec:instab-parabolic-overA}

In this section, let $G$ be a $D$-standard (see \S
\ref{sec:reductive}) reductive group scheme over $\A$.

\subsection{Equidimensional nilpotent sections}

Let $X \in \glie(\A)$ be a section of the Lie algebra $\glie =
\Lie(G)$, and let $C_G(X) = \Stab_G(X)$ be the centralizer of this
section; cf. \S\ref{sub:smoothness}. On base change, the group
$C_G(X)_t$ is the centralizer of $X(t)$ in the algebraic group $G_t$
for each point $t$ of $\Spec(\A)$; according to
\ref{stmt:smooth-centralizers}, each group $C_G(X)_t$ is smooth over
$k(t)$. In general, however, the group scheme $C_G(X)$ will not be
smooth -- or even flat -- over $\A$.

We say that $X$ is \emph{nilpotent} if the value of $X$ at the generic
point $\eta \in \Spec(\A)$ is nilpotent -- i.e. if $X(\eta) \in
\glie(K)$ is nilpotent.

\begin{stmt}
  If $X$ is nilpotent, then also the value $X(t) \in \glie(k(t))$ is
  nilpotent for each point $t \in \Spec(\A)$.
\end{stmt}

\begin{proof}
  Let $\lambda$ denote the regular representation of $G$ on $\A[G]$.
  Since $G$ is reductive, it is by definition smooth -- and in
  particular, flat -- over $\A$. Since the coordinate algebra $\A[G]$
  is a flat $\A$-module, we may regard
  $\A[G]$ as a subring of $K[G]$. Since $X$ is nilpotent, the operator
  $a(\eta)=d\lambda(X(\eta)):K[G] \to K[G]$ is locally nilpotent; i.e. for
  each $f \in K[G]$, we have $a(\eta)^{N(f)}f = 0$ for some $N(f) > 0$.

  Let $t \in \Spec(\A)$, and consider the localization $\A_t = \A_\pp$
  where $\pp \subset \A$ is the prime ideal that ``is'' the point $t$.
  Then we have $\A[G] \subset \A_t[G] \subset K[G]$, and $a(\eta)$
  restricts to a locally nilpotent endomorphism $a$ of $\A[G]$ and of
  $\A_t[G]$. Since $k(t)$ is a quotient of $\A_t$, it follows that
  $a(t) = d\lambda(X(t)):k(t)[G] \to k(t)[G]$ is locally nilpotent, so
  that $X(t)$ is indeed nilpotent as required.
\end{proof}

We say that a nilpotent section $X \in \glie(\A)$ is
\emph{equidimensional} if $\dim C_G(X)_t$ is constant for each $t \in
\Spec(\A)$. For example, if $G = \GL_3$ and $\A$ is a discrete
valuation ring with uniformizing element $\pi$, consider the nilpotent
sections
\begin{equation*}
  X_1=\begin{pmatrix}
    0 & 1 & 0 \\
    0 & 0 & 1 \\
    0 & 0 & 0
  \end{pmatrix},
  X_2=\begin{pmatrix}
    0 & 1 & 0 \\
    0 & 0 & \pi \\
    0 & 0 & 0
  \end{pmatrix} \in \mathfrak{gl}_3(\A).
\end{equation*}
Then $X_1$ is equidimensional, while $X_2$ is not.

\begin{stmt}
  Let $\eta,s \in \Spec(\A)$ be respectively the generic point and the
  closed point.  If 
  \begin{equation*}
    \dim C_G(X)_\eta = \dim C_G(X)_s,    
  \end{equation*}
  then $X$ is equidimensional.
\end{stmt}

\begin{proof}
  Since $C_G(X)$ is the fiber product $G \times_\glie
  \Spec(\A)$, it is a scheme of 
  finite type over $\A$. The assertion now follows from
  \ref{stmt:dim-behaviour}.
\end{proof}

\subsection{Smoothness}
\label{sub:smooth-over-A}

If $\LL$ is a free $\A$-module of finite rank $d$, we can regard $\LL$
as an $\A$-scheme isomorphic to $\Aff^d$. Moreover, we may consider
the $\A$-scheme $\Proj(\LL)$ given for each commutative $\A$-algebra
$\Lambda$ by
\begin{equation*}
  \Proj(\LL)(\Lambda) = \text{set of those $\Lambda$-direct summands of 
    $\LL \tensor_{\A} \Lambda$ having rank 1.}
\end{equation*}
Then $\Proj(\LL)$ is  isomorphic to $\Proj^{d-1}$.  For each $s
\in \Spec(\A)$, the scheme $\Proj(\LL)_s$ obtained by
base-change is just the projective space of the $k(s)$-vector space
$\LL_s = \LL \tensor_{\A} k(s)$.

We are going to consider the $\A$-schemes $\glie$ and $\Proj(\glie)$
where $\glie=\Lie(G)$. Of course, $G$ acts on $\glie$ by the adjoint
representation. If $Y \in \glie(A)$, we write $C_G(Y)$ for the
stabilizer $\Stab_G(Y)$ of the sections $Y$.

The adjoint action of $G$ on $\glie$ determines also an action of $G$
on the projective space $\Proj(\glie)$. If $Y \in \glie(\A)$ is a
section whose image in $\glie(k)$ is non-zero, where $k$ is the
residue field of $\A$, then $\A Y$ is an $\A$-direct summand of
$\glie(\A)$, so it determines a section $[Y] \in \Proj(\glie)(\A)$. We
write $N_G(Y) = \Stab_G([Y])$ for this stabilizer.

Let $X \in \glie(\A)$ be a non-zero equidimensional nilpotent section.
Since $X$ is equidimensional, evidently $X(s) \ne 0$; thus $X$ 
determines a section $[X] \in \Proj(\glie)(\A)$.
\begin{prop}
  The subgroup schemes $C_G(X)$ and $N_G(X)$ of $G$ are smooth over
  $\A$.
\end{prop}

\begin{proof}
  For each $t \in \Spec(\A)$, we know that the dimension of $N_G(X)_t$
  is one more than the dimension of $C_G(X)_t$; cf.
  \cite{jantzen-nil}*{\S 5.3}.  We know from
  \ref{stmt:smooth-centralizers} that $C_G(X)_t$ is smooth, and it
  follows from \cite{mcninch-rat}*{Lem. 23} that $N_G(X)_t$ is
  smooth.  Thus the Proposition follows from \ref{stmt:smoothness-II}
  using equidimensionality.
\end{proof}

\subsection{Richardson sections}
Let $X = \Aff^n_{/\Z}$ be affine $n$-space over $\Z$ for some $n \ge
1$, and let $S = \{p_1,\dots,p_n\}$ be a finite set of $n$ distinct
prime numbers. We regard the $p_i$ as points of $\Spec(\Z)$, and we
write $\xi$ for the generic point of $\Spec(\Z)$.

Suppose that we are given an open $\Q$-subscheme $U_0$ of the generic
fiber $X_\xi$, and that for each $p \in S$, we are given an open
$\FF_p$-subscheme $U_p$ of the fiber $X_p$.

For a regular function $f \in \Z[X]$ and a field $E$, we write $f_E$
for the corresponding regular function $f \tensor 1$ in $E[X] = \Z[X]
\tensor_\Z E$.

\begin{stmt}
  \label{stmt:regular-function-open-on-fibers}
  There is a regular function $f \in \Z[X]$ such
  that 
  \begin{enumerate}
  \item[(i)] the distinguished open subset of $X_\xi$ determined by
    the non-vanishing of $f_\Q$ is contained in $U_0$, and
  \item[(ii)] for $p \in S$, the distinguished open subset of
    $X_p$ determined by the non-vanishing of $f_{\FF_p}$ is
    contained in $U_p$.
  \end{enumerate}
\end{stmt}

\begin{proof}
  Let $g \in \Q[X]$ be a regular function such that the distinguished
  open subset $D(g)$ of $X_\xi$ determined by the non-vanishing of $g$
  lies in $U_\xi$.  We may evidently replace $g$ by a non-zero integer
  multiple without changing $D(g)$; since $\Z[X] =
  \Z[T_1,\dots,T_n]$ is a factorial domain, we may suppose that $g \in
  \Z[X]$ and that the image of $g$ in $\FF_p[X]$ is non-zero for
  each prime $p$.

  For $p \in S$, let $h_p \in \FF_p[X]$ be a non-zero regular function
  such that the distinguished open subset $D(h_p)$ of $X_p$ determined
  by $h_p$ lies in the open subscheme $U_p$.  Since the natural
  mapping $\Z[X] \to \prod_{p \in S} \FF_p[X]$ is surjective by the
  Chinese Remainder Theorem, we may find $h \in \Z[X]$ whose image in
  $\FF_p[X]$ is $h_p$ for each $p \in S$.

  Now put $f = g \cdot h \in \Z[X]$. Then for each $p \in S$, the
  image $f_{\FF_p}$ is non-zero; since $h_p \mid f_{\FF_p}$, the
  distinguished open subset of $X_p$ determined by the non-vanishing
  of $f_{\FF_p}$ is contained in $U_p$.  Moreover, since $g \mid f =
  f_{\Q}$, the distinguished open subset of $X_\xi$ determined by the
  non-vanishing of $f_{\Q}$ is contained in $U_0$, as required.
\end{proof}

Let now $\A$ be a local, normal, Noetherian domain, and suppose that
$G$ is a split reductive group over $\A$, with split maximal torus
$T$.
\begin{stmt}
  \label{stmt:first-Richardson-section}
  Let $P \subset G$ be a parabolic subgroup scheme containing $T$. 
  \begin{enumerate}
  \item[(a)] There is a regular function $f \in \A[\Lie(R_uP)]$ such
    that for each $t \in \Spec(\A)$, the distinguished open subset of
    the $k(t)$-scheme $\Lie(R_uP)_t$ determined by the non-vanishing
    of $f_{k(t)}$ is contained in the Richardson orbit of $P_t$ on
    $\Lie(R_uP)_t$.
  \item[(b)] If the residue field of $\A$ is infinite, there is a
    section $X \in \Lie(R_uP)(\A)$ such that $X(t)$ is a Richardson
    element for $P_t$ for each $t \in \Spec(\A)$.
  \end{enumerate}
\end{stmt}

\begin{proof}
  There is a split reductive group scheme $G_0$ over $\Z$ and a split
  maximal torus $T_0$ such that $G = G_{0/\A}$ and $T=T_{0/\A}$.
  
  Now choose an $A$-homomorphism $\phi:\Gm \to T$ such that $P =
  P(\phi)$ as in \ref{stmt:P-is-P(phi)}. Since $T_0$ is a split torus
  over $\Z$, there is a $\Z$-homomorphism $\psi:\Gm \to T_0$ such that
  $\phi = \psi_{/\A}$. Writing $P_0$ for the parabolic $\Z$-subgroup
  scheme of $G_0$ determined by $\psi$, we have $P = P_{0/\A}$.

  On the geometric fibers, it follows from the finiteness of the
  number of nilpotent $G_{0,\bar p}$-orbits \footnote{That finiteness
    is true in all characteristics, though the proof in bad
    characteristic is ``case-by-case'' at present. See
    \cite{jantzen-nil}*{\S 2.8}.} in $\Lie(G_0)_{\bar p}$ that
  $P_{0,\bar p}$ has an open orbit (the Richardson orbit) in
  $\Lie(R_u(P_0))_{\bar p}$ for each $p \in \Spec(\Z)$. Using
  \cite{springer-LAG}*{Prop.  11.2.8} one knows that this open orbit
  is obtained by base change from an open $k(p)$-subscheme for each $p
  \in \Spec(\Z)$.

  Let $p$ denote the characteristic of the residue field of $\A$; if
  $p>0$, let $S = \{p\}$; otherwise, let $S = \emptyset$.  Now use
  \ref{stmt:regular-function-open-on-fibers} applied to $X =
  \Lie(R_uP_0)$ and the set $S$ to find a regular function $f \in
  \Z[\Lie(R_uP_0)]$ whose image in $\A[\Lie(R_uP)]$ has the required
  properties. This proves (a).

  For (b), since $k(s) = k$ is infinite, where $s \in \Spec(\A)$ is
  the closed point, one may choose an element $Y \in \Lie(R_uP)(k)$
  such that the regular function $f_{k}$ does not vanish at $Y$. Let
  $X \in \Lie(R_uP)(\A)$ be any section such that $X(s) = Y$. Then
  evidently the value $f(X)$ is a unit in $\A$; it now follows from
  (a) that $X(t)$ is a Richardson element for $P_t$ for each $t \in
  \Spec(\A)$, as required.
\end{proof}

\subsection{Existence of equidimensional nilpotent sections}
\label{sub:existence-equidim}

Assume throughout this section that $\A$ is a normal, local,
Noetherian domain with \emph{infinite residue field} $k$, and that the
reductive group scheme $G$ is \emph{split} over $\A$, with split
maximal torus $T \subset G$. We suppose that $G$ is
\emph{$D$-standard.}

Let $L \subset G$ be a Levi factor of some parabolic subgroup scheme
of $G$, and suppose that $T \subset L$. We remark that $L$ itself is
$D$-standard.

Since $T$ is a split torus, we may choose an isomorphism $T \simeq
D_\A(X)$ where $X=X(T)$ is the free Abelian group $\Z^r$; then $X(T)$
identifies with the group of characters $\Hom_\A(T,\Gm)$.  From the
roots of $L$ with respect to $T$, choose a system of positive roots
$R^+ \subset X$ and a basis of the roots $\Pi \subset R^+$.

Now write $\Der(L)=L'$ for the derived subgroup scheme as in
\ref{sub:derived-group}, and write $T'$ for the maximal split torus of
$L'$ contained in $T$; cf.  \ref{stmt:max-torus-of-Der(G)}.
Then $X(T')$ contains $\Z R$ as a subgroup of finite index.

\begin{stmt}
  \label{stmt:Q-L-cochar}
  Let $Q_0 \subset L_t$ be a distinguished parabolic subgroup
  containing $T_t$ for some $t \in \Spec(\A)$.  
  \begin{enumerate}
  \item[(a)] There is a parabolic subgroup scheme $Q \subset L$ such that
    $Q_0 = Q_t$.
  \item[(b)] $Q_x$ is a distinguished parabolic subgroup of $L_x$ for every
    $x \in \Spec(\A)$.
  \item[(c)] Let $I \subset \Pi$ be defined by the condition $\alpha \in I
    \iff \Lie(Q_0)_{-\alpha} \ne 0$. Then there is a unique
    $\A$-homomorphism $\phi:\Gm \to T$ such that $\langle \alpha, \phi
    \rangle = 2$ for  $\alpha \in \Pi \setminus I$
    and $\langle \alpha,\phi \rangle = 0$ for $\alpha \in I$. Moreover,
    $Q$ is the parabolic subgroup of $L$ determined by $\phi$.
  \end{enumerate}
\end{stmt}

\begin{proof}
  Since $L$ is $D$-standard, the characteristic of $k(x)$ is good
  for the derived group of $L_x$ for every $x \in \Spec(\A)$.  It then
  follows from \cite{jantzen-nil}*{Lem. 5.2} that the homomorphism
  $\Z R \to \Z$ given by the rule in (c) determines an
  $\A$-homomorphism $\phi:\Gm \to T'$.

  Let $Q$ be the parabolic subgroup of $L$ determined by $\phi$.  (a)
  is then clear, and (b) follows from \cite{jantzen-nil}*{\S4.10(2)}.
\end{proof}

\begin{stmt}
  \label{stmt:Richardson-section}
  Let $\phi:\Gm \to T$ be the cocharacter of \ref{stmt:Q-L-cochar}(c).
  There is a section
  \begin{equation*}
    X \in \Lie(L)(\phi;2)(\A)
  \end{equation*}
  such that
  \begin{enumerate}
  \item[(a)] $X(t)$ is a Richardson element for $Q_t$,
  \item[(b)] $\phi_t$ is associated with $X(t)$, and
  \item[(c)] the Bala-Carter datum of $X(t)$ is $(L_t,Q_t)$.
  \end{enumerate}
  for each $t \in \Spec(\A)$. Moreover, $X$ is an
  equidimensional nilpotent section of $\Lie(G)$.
\end{stmt}

\begin{proof}
  Since the residue field of $\A$ is assumed to be infinite, we may use
  \ref{stmt:first-Richardson-section} to find a section $Y \in
  \Lie(R_uQ)(\A)$ such that $Y(t)$ is a Richardson element for $Q_t$
  for each $t \in \Spec(\A)$.

  It follows from \cite{jantzen-nil}*{Lem. 5.2 and Lem. 5.3} that
  $\phi_t$ is almost associated with $X(t)$ for each $t \in
  \Spec(\A)$.

  Since $Q$ is the parabolic subgroup determined by $\phi$, we know
  that
  \begin{equation*}
    \Lie(R_uQ) = \sum_{i \ge 1} \Lie(L)(\phi;i).
  \end{equation*}
  Thus, we may write $Y = \sum_{i \ge 1} Y_i$ with $Y_i \in
  \Lie(L)(\phi:i)(\A)$.

  It follows from \ref{stmt:dense-P-orbit} that for each $t \in
  \Spec(\A)$, the element $Y_2(t)$ is Richardson for $Q_t$, and the
  cocharacter $\phi_t$ is associated with $Y_2(t)$; in particular, if
  we set $X = Y_2$ then (a), (b), and (c) hold for $X$.

  Write $P$ for the parabolic subgroup scheme $P(\phi) \subset G$.
  Since $\phi_t$ is associated with $X(t)$ for each $t \in \Spec(\A)$,
  we know  $P_t$ to be the instability parabolic of $X(t)$, so that
  -- by \ref{stmt:assoc-cochar} -- we have $C_G(X)_t \subset P_t$ for
  each point $t$ of $\Spec(\A)$.  Now, the $P_t$-orbit of $X(t)$ is
  dense in $\sum_{j \ge 2} \glie(\phi;j)_t$ by
  \ref{stmt:dense-P-orbit}. It follows that the centralizer of $X$ in
  $P$ has constant dimension on the fibers of $\Spec(\A)$, so that $X$
  is indeed equidimensional.  
\end{proof}

\begin{theorem}
  Let $t \in \Spec(\A)$ and let $Y \in \glie(k(\bar{t}))$ be a
  nilpotent element. Then there is a Levi subgroup scheme $L$ of a
  parabolic subgroup scheme of $G$, an $\A$-homomorphism
  $\Phi:\Gm \to L$, and a nilpotent section $X \in \Lie(L)(\Phi;2)(\A)$ for
  which the following conditions hold:
  \begin{enumerate}
  \item[(a)] $X$ is an equidimensional nilpotent section of $\glie$,
  \item[(b)] $X(t)$ is $G_{\bar{t}}$-conjugate to Y,
  \item[(c)] for each $u \in \Spec(\A)$, the Bala-Carter datum of
    $X(\bar{u})$ is $(L_{\bar u},Q_{\bar u})$, where $Q$ is the
    parabolic subgroup scheme $P_L(\Phi)$ of $L$ determined by $\Phi$.
  \item[(d)] $P_u$ is the instability parabolic of $G_u$ determined by
    $X(u)$ for each $u \in \Spec(\A)$, where $P = P_G(\Phi)$ is the
    parabolic subgroup scheme of $G$ determined by $\Phi$.
  \end{enumerate}
  In particular, $\Phi_u$ is a cocharacter of $G_u$ associated with
  $X(u)$ for each $u \in \Spec(\A)$.
\end{theorem}

\begin{proof}
  Recall that $T$ is a fixed split maximal torus of $G$.  Suppose that
  $(L_0,Q_0)$ is the Bala-Carter datum of $Y$; thus $L_0$ is a Levi
  subgroup of a parabolic subgroup of $G_t$, $Y$ is distinguished in
  $\Lie(L_0)$, and $Q_0$ is a distinguished parabolic subgroup of
  $L_0$. Since we work up to geometric conjugacy, we may as well
  suppose that $L_0$ is defined over $k(t)$, that $L_0$ contains
  $T_t$, and that $T_t$ contains a maximal torus of the centralizer in
  $G$ of $Y$.  This last condition shows that $L_0$ is the centralizer
  of the image of some cocharacter of $T_t$. Since $T$ is a split
  torus, this cocharacter arises by base change from an
  $\A$-homomorphism $\Phi:\Gm \to T$; in view of
  \ref{stmt:Levi-of-P-via-Phi}, there is a Levi subgroup scheme $L$ of
  a parabolic subgroup scheme of $G$ for which $L_t = L_0$.  Now use
  \ref{stmt:Q-L-cochar} to see that $Q=P_L(\Phi)$ is a distinguished
  parabolic subgroup scheme of $L$ for which $Q_t$ is $L_{\bar
    \eta}$ conjugate to $Q_0$. Finally, use
  \ref{stmt:Richardson-section} to find an equidimensional nilpotent
  section $X \in \Lie(L)(\Phi;2)(\A)$ for which $X(u)$ has Bala-Carter datum
  $(L_u,Q_u)$ for each $u \in \Spec(\A)$.  Then $X(t)$ is $G_{\bar
    t}$-conjugate to $Y$. Thus (a), (b) and (c) hold.
  
  By \ref{stmt:Richardson-section}, $\Phi_u$ is a cocharacter of $G_u$
  associated with $X(u)$ for each $u \in \Spec(\A)$. Denoting by $P$
  the parabolic subgroup scheme $P_G(\phi)$, we conclude that $P_u$ is
  the instability parabolic of $X(u)$ for each $u \in \Spec(\A)$; thus
  (d) holds as well.
\end{proof}

\subsection{The instability parabolic of $X$}
\label{sub:instability-parabolic-over-A}
Let $X \in \glie(\A)$ be an equidimensional nilpotent section. Let
$\eta \in \Spec(\A)$ be the generic point, and let $P_0 \subset
G_\eta$ by the instability parabolic subgroup determined by $X(\eta)$.

\begin{prop}
  There is a unique $\A$-parabolic subgroup scheme $P \subset G$ such
  that $P_0=P_\eta$.
\end{prop}

\begin{proof}
  Unicity follows from \ref{stmt:parabolic-determined-by-P_eta}(a).
  For existence, first suppose that $\A$ is a discrete valuation ring.
  In that case the conclusion of the Proposition is a consequence of
  \ref{stmt:parabolic-over-dvr}.

  Since the scheme $\Par$ of parabolic subgroups of $G$ is projective
  \ref{stmt:parabolic-scheme-projective}, the Proposition now follows
  in the general case from
  \ref{stmt:section-of-projective-over-normal-A}.
\end{proof}

\begin{rem}
  The conclusion of the Proposition has already been observed for the
  nilpotent sections obtained using Theorem
  \ref{sub:existence-equidim}.
\end{rem}

\subsection{\'Etale local existence of associated cocharacters over $\A$}
\label{sub:existence-of-assoc-cochar-over-A}

With notation as before, write $P \subset G$ be the parabolic subgroup
scheme for which $P_\eta$ is the instability parabolic of $X(\eta)$.
\begin{stmt}
  \label{stmt:dense-Pk-orbit} 
  \begin{enumerate}
  \item[(i)] There is an $\A$-homomorphism $\Phi:\Gm \to P$ such that
    the cocharacter $\Phi_\eta$ is almost-associated with $X(\eta)$.
  \item[(ii)] For each $t \in \Spec(\A)$, the $P_t$-orbit of $X(t)$ is
    separable and dense in $\sum_{i \ge 2} \glie(\Phi;i)_t$.
  \end{enumerate}
\end{stmt}

\begin{proof}
  Using \ref{stmt:parabolic-has-levi}, we choose a maximal torus $T
  \subset P$. Let $\Phi_0$ be the unique cocharacter of $T_\eta$ which
  is almost-associated to $X(\eta)$.  It follows from
  \ref{stmt:cochar-over-A} that there is an $\A$-homomorphism
  $\Phi:\Gm \to T$ inducing $\Phi_0$ on base-change; this proves (i).

  We now prove (ii). Since $\Phi_\eta$ is almost associated with
  $X(\eta)$, the $P_\eta$-orbit of $X(\eta)$ is dense in $\sum_{i \ge
    2} \glie(\Psi;j)_\eta$ by \ref{stmt:dense-P-orbit}.  In
  particular, $X$ may be regarded as an $\A$-section of $\sum_{i \ge
    2} \glie(\Psi;j)$ and so $X(t)$ is a section of $\sum_{i \ge 2}
  \glie(\Psi;j)_t$.
  
  Write $d$ for the $\A$-rank of the free $\A$-module $\sum_{i \ge 2}
  \glie(\Psi;j)$. Since the centralizer of $X(\eta)$ in $G_\eta$ is contained
  in $P_\eta$, we have by assumption that
  \begin{equation*}
    \dim P_\eta - \dim C_G(X)_\eta = d.
  \end{equation*}

  Now, we certainly have $\dim C_P(X)_t \le \dim C_G(X)_t$.
  Since the $P_t$-orbit of $X(t)$ lies in $\sum_{i \ge 2}
  \glie(\Psi;j)_t$, and since $X$ is equidimensional, this orbit has
  dimension
  \begin{equation*}
    \dim P_t - \dim C_P(X)_t \ge 
    \dim P_t - \dim C_G(X)_t = d.
  \end{equation*}
  For dimension reasons, we conclude that the $P_t$-orbit of $X(t)$ is
  dense in $\sum_{i\ge 2} \glie(\Psi;i)_t$.  Since $C_G(X)_t$ is
  smooth \ref{stmt:smooth-centralizers}, the dimension of the
  centralizer of $X(t)$ in the Lie algebra $\glie_t$ coincides with
  $\dim C_G(X)_t$.  It follows that the dimension of the
  centralizer of $X(t)$ in the Lie algebra $\Lie(P)_t$ must coincide
  with $\dim C_P(X)_t$, so that the $P_t$-orbit of $X(t)$ is
  indeed separable.
\end{proof}

Let $M$ be a free $\A$-module of finite rank, and write $M_\eta = M
\tensor_\A K$.
\begin{stmt}
  \label{stmt:equal-if-so-on-generic-fiber}
  If $N,N' \subset M$ are $\A$-direct summands of $M$, then $N =
  N'$ if and only if $N_\eta = N'_\eta$.
\end{stmt}

\begin{proof}
  Indeed, for any $\A$-direct summand $L$ of $M$, we have $L = L_\eta
  \cap M$.  Thus $N_\eta = N'_\eta$ indeed implies that $N = N'$; the other
  implication is even simpler.
\end{proof}

\begin{prop}
  Let $P_1 \subset G$ be the $\A$-parabolic
  subgroup scheme for which $P_{1,\eta}$ is the instability parabolic of
  $X(\eta)$.  Let $\Psi:\Gm \to P_1$ be an $\A$-homomorphism such that
  $\Psi_\eta$ is almost associated to $X(\eta)$. Then 
  \begin{enumerate}
  \item[(i)] $P_{1,t}$ is the instability parabolic of $X(t)$ for each
    $t \in \Spec(\A)$
  \item[(ii)] $\Psi_t$ is almost associated to $X(t)$ for each $t \in
    \Spec(\A)$.
  \item[(iii)] There is a finite, \'etale, local extension $\B \supset
    \A$ and a section $g \in P_1(\B)$ such that if we put
    \begin{equation*}
      \Phi=\Int(g) \circ \Psi:\Gm    \to P_1,
    \end{equation*}
    then $\Phi_t$ is a cocharacter associated to $X(t)$ for each $t
    \in \Spec(\A)$.
  \end{enumerate}
\end{prop}

\begin{proof}
  Since $P_{1,\eta} = P_1(\Psi)_{/\eta}$, it follows from the
  uniqueness assertion in Proposition
  \ref{sub:instability-parabolic-over-A} that $P_1 = P(\Psi)$.  Also
  notice that (i) and (ii) are consequences of (iii); we will just
  prove (iii).

  The statement is unchanged if we replace $\A$ by a finite, \'etale,
  local extension; thus, we may and will suppose that $G$ is split,
  say with split maximal torus $T$.

  Using Theorem \ref{sub:existence-equidim}, we locate a Levi subgroup
  scheme $L$ of a parabolic subgroup scheme of $G$, an
  $\A$-homomorphism $\Phi:\Gm \to L$ and a section $Y \in
  \Lie(L)(\Phi;2)(\A)$ such that 
  \begin{itemize}
  \item $Y$ is an equidimensional nilpotent section of $\Lie(G)$,
  \item $Y(\eta)$ is $G_{\bar{\eta}}$-conjugate to $X(\eta)$,
  \item $(L_u,Q_u)$ is the Bala-Carter datum of $Y(u)$ for each $u \in
    \Spec(\A)$, where $Q = P_L(\Phi)$ is the parabolic subgroup scheme
    of $L$ determined by $\Phi$, and
  \item for every $u \in \Spec(\A)$, $P_u$ is the parabolic subgroup
    of $G_u$ associated with $Y(u)$ and $\Phi_u$ is a cocharacter of
    $G_u$ associated with $Y(u)$, where $P = P_G(\Phi)$ is the
    parabolic subgroup scheme of $G$ determined by $\Phi$.
  \end{itemize}
  We may evidently suppose that $Q$ contains the split maximal torus
  $T$ of $G$.

  We know that $P_{1,\eta}$ and $P_\eta$ are $G_{\bar
    \eta}$-conjugate.  After possibly replacing $\B$ by a finite
  \'etale, local extension, we may suppose that $P_{1,\eta}$ and
  $P_\eta$ are conjugate by an element in $G(k(\eta))$.  Using
  \cite{SGA3}*{Exp. XXVI, Cor. 5.5 (i)} we see that $P_1$ and $P$ are
  $G(\A)$-conjugate; thus we may and will suppose that $P_1=P$. Then
  both cocharacters $\Phi_\eta$ and $\Psi_\eta$ are almost associated
  with $X(\eta)$.

  It follows that $\sum_{i \ge 2} \glie(\Psi;i)_\eta = \sum_{i \ge 2}
  \glie(\Phi;i)_\eta$; using \ref{stmt:equal-if-so-on-generic-fiber},
  we can now conclude that $\sum_{i \ge 2} \glie(\Psi;i) = \sum_{i \ge
    2} \glie(\Phi;i).$

  By \ref{stmt:dense-Pk-orbit} the $P_t$-orbits of $X(t)$ and of
  $Y(t)$ are separable and dense in 
  \begin{equation*}
   \sum_{j \ge 2} \glie(\Phi;j)_t =
   \sum_{j \ge 2} \glie(\Psi;j)_t  
  \end{equation*}
  for each point $t$ of $\Spec(\A)$.  
  Using \ref{stmt:transport-via-etale-ext} we may find a finite,
  \'etale, local extension $\B$ of $\A$ such that $X$ and $Y$ are
  conjugate by an element of $P(\B)$; we may and will replace $\A$
  with $\B$ so that $X$ and $Y$ are conjugate by an element of
  $P(\A)$. Thus, we may and will suppose that $X=Y$.

  Since the centralizers of $\Phi$ and of $\Psi$ are Levi subgroup
  schemes of $P$ \ref{stmt:Levi-of-P-via-Phi}, we may find maximal
  tori $T_1,T_2 \subset P$ such that $\Phi$ factors through the
  inclusion of $T_1$ in $P$ and such that $\Psi$ factors through the
  inclusion of $T_2$ in $P$. Since $T_1$ and $T_2$ are locally
  conjugate for the \'etale topology of $P$ \cite{SGA3}*{Exp. XII,
   Thm 1.7}, after replacing $\A$ by a finite, \'etale, local
  extension, the maximal tori $T_1$ and $T_2$ are conjugate by an
  element of $P(\A)$. Thus we may suppose that $T_1 = T_2$; but then
  $\Phi_\eta = \Psi_\eta$ by \ref{stmt:almost-assoc-exists}.  It now
  follows that $\Phi=\Psi$.  But then one knows for each $t \in
  \Spec(\A)$ that $\Phi_t=\Psi_t$ is associated with $Y(t) = X(t)$ and
  the proof is complete.
\end{proof}

\subsection{Maximal tori and Levi factors.}
\label{sub:max-tori-and-levi}
We are going to prove in this section the main Theorem regarding the
existence of a Levi factor of the centralizer of an equidimensional
nilpotent section. We first require a preliminary observation.

Let $H$ be a smooth group scheme over $\A$. For $t \in \Spec(\A)$, let
$\rho_r(t) = \rho_{r,H}(t)$ be the dimension of a maximal torus of the
$k(\bar{t})$-group $H_{\bar t}$ for some (hence any) geometric point
$\bar{t}$ above $t$.

\begin{stmt}
  \label{stmt:torus-condition}
  The following are equivalent:
  \begin{enumerate}
  \item[(a)] The function $\rho_r$ is constant on $\Spec(\A)$.
  \item[(b)] $\rho_r(s) = \rho_r(\eta)$ where $s$ and $\eta$ are
    respectively the closed point and the generic point of $\Spec(\A)$.
  \item[(c)] Locally in the \'etale topology, $H$ has a maximal torus.
  \end{enumerate}
\end{stmt}

\begin{proof}
  The equivalence of (a) and (c) follows from \cite{SGA3}*{Exp. XII,
    Thm. 1.7(b)}, while the equivalence of (a) and (b) follows
  from the lower semi-continuity of $\rho_r$ on $\Spec(\A)$; cf.
  \emph{loc. cit.}  Thm. 1.7(a).
\end{proof}

\begin{theorem}
  Let $G$ be a $D$-standard reductive group scheme over $\A$, let
  $X \in \glie(\A)$ be an equidimensional nilpotent section, let $C =
  C_G(X)$, and let $N = N_G(X)$.  
  \begin{enumerate}
  \item[(a)] There is a finite, \'etale, local extension $\B \supset
    \A$ and a $\B$-homomorphism $\phi:\Gm \to G_{/\B}$ such that the
    cocharacter $\phi_t$ of $G_t$ is associated to $X(t)$ for
    each $t \in \Spec(\B)$.
  \item[(b)] $C$ has a Levi factor locally in the \'etale topology of
    $\Spec(\A)$.
  \item[(c)] $C$ has a maximal torus locally in the \'etale topology of $\Spec(\A)$.
  \item[(d)] $N$ has a maximal torus locally in the \'etale topology of $\Spec(\A)$.
  \end{enumerate}
\end{theorem}

\begin{proof}
  Let $P$ be the $\A$-parabolic subgroup scheme of $G$ for which
  $P_\eta$ is the instability parabolic for $X(\eta)$; see Proposition
  \ref{sub:instability-parabolic-over-A}.  Now, (a) has been proved
  already in Proposition \ref{sub:existence-of-assoc-cochar-over-A}.
  

  In order to prove (b), (c) and (d), we may and will replace $\A$ by
  a finite, \'etale, local extension; thus we may suppose by part (a)
  that $\phi:\Gm \to G$ is an $\A$-homomorphism for which $\phi_t$ is
  associated with $X(t)$ for all $t \in \Spec(\A)$. 

  The centralizer $L$ of the image of $\phi$ in $C$ is a (closed)
  subgroup scheme of $C$, and $L$ is smooth over $\B$; cf.
  \ref{stmt:mult-centralizer}.  Moreover, it follows from
  \ref{stmt:assoc-cochar}(d) that $L_t$ is a Levi factor of $C_t$ for
  each $t \in \Spec(\A)$, whence (b).

  Since the subgroup scheme $L^o$ is reductive, one knows that $L$ --
  and hence $C$ -- has a maximal torus by
  \ref{stmt:reductive-has-max-torus}; this proves (c).

  According to \cite{jantzen-nil}*{\S 5.3}, one knows for each $t \in
  \Spec(\A)$ that $N_G(X)_{\bar t}$ is the product of $C_G(X)_{\bar
    t}$ with the image of any cocharacter of $G_{\bar t}$ associated
  with $X$. Since the image of such a cocharacter centralizes some
  maximal torus in $C_G(X)_{\bar t}$, it follows that $\rho_{r,C}(t) +
  1 = \rho_{r,N}(t)$.  Using \ref{stmt:torus-condition} it is now
  clear that (d) is a consequence of (c).
  
\end{proof}

For later use, we observe that the proof of part (b) of the preceding
Theorem actually proves the first assertion of the following:
\begin{stmt}
  \label{stmt:phi=>Levi}
  Assume that there is an $\A$-homomorphism $\phi:\Gm \to G$ such
  that the cocharacter $\phi_t$ of $G_t$ is associated to $X(t)$ for
  each $t \in \Spec(\A)$.
  \begin{enumerate}
  \item the centralizer of the image of $\phi$ in $C$ is a Levi factor
    $L$ of $C$.
  \item There is a smooth retraction $\rho:C \to L$ in the sense of
    \ref{stmt:levi-extra}; in particular, writing $R = \ker \rho$,
    there is an isomorphism of $\A$-group schemes $C \simeq L \ltimes
    R$.
  \end{enumerate}
\end{stmt}

\begin{proof}
  We have  observed that (1) was proved already. For the second
  assertion, write $\A[C] = \bigoplus_{n \in \Z} \A[C]_n$ as a direct
  sum of weight spaces for the action of $\Gm$ on $\A[C]$ given by
  $\Int^* \circ \phi$, where $\Int$ is the action of $C$ on itself by
  inner automorphisms.

  Since $\phi_t$ is associated to $X(t)$, one knows by
  \ref{stmt:assoc-cochar} that $\A[C]_{n,t} = 0$ for any $t \in
  \Spec(\A)$ whenever $n > 0$. It follows that $\A[C]_n = 0$ whenever
  $n > 0$; i.e. 
  \begin{equation*}
    \A[C] = \bigoplus_{n \le 0} \A[C]_n.
  \end{equation*}
  Write $\A[C]_{<0} = \sum_{n < 0} \A[C]_n$.  Then $\A[L] =
  \A[C]/\A[C]_{<0}$, and the inclusion mapping $i:L \to C$ is given by
  the natural surjection $i^*:\A[C] \to \A[C]/\A[C]_{>0}$. The Hopf
  algebra $\A[L]$ identifies naturally with $\A[C]_0$, and the
  inclusion map $\rho^*:\A[C]_0 \to \A[C]$ defines a retraction
  $\rho:C \to L$. Since $\rho_t$ is evidently smooth for all $t \in
  \Spec(\A)$, and since $C$ and $L$ are both flat over $\A$,
  \cite{SGA1}*{Exp. II, Cor. 2.2} shows that $\rho$ is a smooth mapping. In
  view of \ref{stmt:levi-extra}, this completes the proof of (2).
\end{proof}


\begin{cor} With assumptions as before, we have the following:
  \begin{enumerate}
  \item[(a)] Locally in the \'etale topology there are subgroup
    schemes $Q \subset L \subset G$ such that $L$ is a Levi subgroup
    scheme of a parabolic subgroup scheme of $G$, $Q$ is a a
    distinguished parabolic subgroup scheme of $L$, and 
    $(L_{\bar{t}},Q_{\bar{t}})$ is the Bala-Carter datum of
    $X(t)$ for each $t \in \Spec(\A)$.
  \item[(b)] The root datum of the connected component of a Levi
    factor of $C_G(X)_{\bar{t}}$ is constant for $t \in \Spec(\A)$.
  \end{enumerate}
\end{cor}

\begin{proof}
  For the proof of the corollary, we may replace $\A$ by a finite,
  \'etale, local extension; applying the Theorem for $G$, we may
  suppose that $C=C_G(X)$ has a maximal torus $T$.  Now let
  $L=C_G(T)$; then $L$ is a Levi factor of a parabolic subgroup scheme
  of $G$, and $X \in \Lie(L)(\A)$. Since $T$ is a maximal torus of
  $C_L(X)$, the Theorem applies also to $L$. Thus, we may suppose that
  there is an $\A$-homomorphism $\phi:\Gm \to L$ such that $\phi_t$ is
  a cocharacter of $L_t$ which is associated to $X(t)$ for each $t \in
  \Spec(\A)$.

  For (a), let $Q = P_L(\phi)$ be the parabolic subgroup scheme of $L$
  determined by $\phi$. Since $\phi_t$ is associated with $X(t)$, one
  knows that $Q_t$ is the instability parabolic subgroup of $L_t$
  determined by $X(t)$ \ref{stmt:assoc-cochar}. Since $T_t$ is a
  maximal torus of $C_G(X)_t$ for each $t \in \Spec(\A)$, it is clear
  that $X(t)$ is distinguished in $\Lie(L)_t$. Thus indeed $(L_t,Q_t)$
  is the Bala-Carter datum of $X(t)$.

  For (b), note that $L_t$ is a Levi factor of $C_G(X)_t$ for each $t
  \in \Spec(\A)$. So (b) follows from \cite{SGA3}*{III Exp. XXII Prop.
    2.8}.
\end{proof}

\subsection{Proof of Theorem \ref{theorem:levi-factor}}
\label{sub:proof-levi-factor}

Recall that $E$ is an algebraically closed field of characteristic 0,
and that $F$ is an algebraically closed field of characteristic $p>0$.
Theorem \ref{theorem:levi-factor} is a consequence of the following
more general result:
\begin{theorem}
  Let $G_E$ and $G_F$ be reductive groups respectively over $E$ and
  over $F$, assume that the root datum of $G_E$ coincides with that of
  $G_F$, and assume that $G_F$ is $D$-standard. Let $X_E \in \glie_E$,
  $X_F \in \glie_F$ be nilpotent elements with the same Bala-Carter
  data, let $C_E$ and $C_F$ be their respective centralizers, and let
  $L_E \subset C_E$ and $L_F \subset C_F$ be Levi factors
  \ref{stmt:assoc-cochar}. Then the root datum of $L_E^o$ may be
  identified with that of $L_F^o$.
\end{theorem}

\begin{proof}
  Let $\A$ be the ring of Witt vectors \footnote{The Witt vectors are
    just a convenient choice. In fact, one can use instead any normal,
    local Noetherian domain $\A$ with infinite residue field of
    characteristic $p>0$ whose field of fractions has characteristic
    0.}  \cite{serre-local-fields}*{II \S6} with residue field an
  algebraic closure of the finite field $\mathbf{F}_p$.  Using
  \ref{stmt:choose-favorite-fields} we see that it is enough to prove
  Theorem \ref{theorem:levi-factor} after replacing $F$ by the residue
  field of $\A$ and $E$ by an algebraic closure of the field of fractions
  of $\A$, and after replacing $X_F$ and $X_E$ by nilpotent elements with the given
  Bala-Carter datum.

  Let $G$ be a split reductive group scheme over $\A$ with the same
  root datum as $G_F$ -- for the existence, see e.g. \cite{SGA3}*{Exp.
    XXV, Thm. 1.1}.  Then $G_F$ identifies with the closed fiber $G_s$
  of $G$, and $G_E$ identifies with the generic fiber $G_{\bar \eta}$,
  where $\eta$ is the generic point of $\Spec(\A)$.
  
  Use Theorem \ref{sub:existence-equidim} and the Bala-Carter theorem
  to find an equidimensional nilpotent section $X$ for which $X(s)$ is
  conjugate to $X_F$ and for which $X(\eta)$ is (geometrically)
  conjugate to $X_E$. We may and will replace $X_F$ by $X(s)$ and
  $X_E$ by $X(\eta)$.

%

If $C$ denotes the centralizer
in $G$ of the nilpotent section $X$, it follows from part (b) of
Corollary \ref{sub:max-tori-and-levi} that the root datum of a Levi
factor of $C_G(X)_{\bar t} = C_{G_{\bar t}}(X(t))$ is constant for $t
\in \Spec(\A)$. This yields the desired result.   
\end{proof}


\section{The group of components of a group scheme}
\label{sec:component-group}

Again let $\A$ be a Noetherian, normal, local domain. Our goal in the
section following this one is to investigate the groups $C_t/C^o_t$
where $t \in \Spec(\A)$, where $C$ is the group scheme $C_G(X)$ for an
equidimensional nilpotent section $X$, and where $G$ is assumed to be
a  $T$-standard reductive group scheme over $\A$. We first require
some preliminaries, which we study in this section.

Let $H$ be a smooth and separated group scheme over $\A$.  We are
going to study the sheaf-quotient $H/H^o$, which we now describe.

\subsection{Sheaves} 
\label{sub:sheaves}
If $X$ is an $\A$-scheme, an \'etale covering of $X$ is a family of
\'etale $\A$-morphisms $(U_i \to X)_i$ of finite type, such that $X =
\cup_i U_i$.

Following \cite{milne}*{II \S1} we write $\Spec(\A)_{\et}$ for the
(small) \'etale site of  $\Spec(\A)$.  This means first
of all that the underlying category of $\Spec(\A)_{\et}$ is the
category $\Et{\Spec(\A)}$ of all schemes which are \'etale and of
finite type over $\Spec(\A)$; the morphisms of this category are just
the morphisms of $\A$-schemes. Finally, $\Spec(\A)_{\et}$ is this
category together with its (Grothendieck) topology defined by
\'etale coverings of finite type.

A pre-sheaf of groups on $\Spec(\A)_{\et}$ is a contravariant functor
\begin{equation*}
  \F:\Et{\Spec(\A)} \to \mathbf{Groups};
\end{equation*}
the pre-sheaf $\F$ is a sheaf if the sequence (S) of \cite{milne}*{II
  \S1 p. 49} is exact for all coverings in $\Et{\Spec(\A)}$.

\begin{stmt}[\cite{milne}*{II Cor. 1.7 and Rem. 1.12}]
  \label{stmt:group-scheme-sheaf}
  A group scheme $H$ over $\Spec(\A)$ determines a sheaf of groups
  on $\Spec(\A)_{\et}$ by the rule $U \mapsto \operatorname{Mor}_{\A}(H,U)$.
\end{stmt}

We say that a sheaf $\F$ on $\Spec(\A)_{\et}$ is \emph{representable}
if there is a group scheme $H$ such that $\F$ is isomorphic to the
sheaf obtained from $H$ in \ref{stmt:group-scheme-sheaf}.  We
sometimes abuse notation and write $H$ for the sheaf $\F$.


If $\F$ is a pre-sheaf on $\Spec(\A)_{\et}$ and if $x \in \Spec(\A)$,
the stalk $\F_{\bar{x}}$ is given by
\begin{equation*}
  \F_{\bar{x}} = \lim_{\rightarrow} \F(U)
\end{equation*}
where the limit is taken over all \'etale neighborhoods $U$ of $x$
\cite{milne}*{II \S2}.

Let $\A_{\bar{x}}$ be the strict Henselization of
$\A_x$ \cite{milne}*{I \S4}; thus $\A_{\bar{x}}$ is a Henselian
local domain containing $\A_x$, and the residue field of the local
ring $\A_{\bar{x}}$ is the separably closed field $k(\bar{x})$.

\begin{stmt}[\cite{milne}*{II Rem. 2.9(d)}]
  \label{stmt:group-scheme-stalk}
  Let $H$ be an $\A$-group scheme, and let $x \in \Spec(\A)$.  Then
  the stalk $H_{\bar{x}}$ of the sheaf $H$ identifies with the group
  of points $H(\A_{\bar{x}})$.
\end{stmt}

Now let $J \subset H$ be a closed and normal subgroup scheme, and
suppose that $H^o \subset J$. We write $H/J$ for the sheaf on
$\Spec(\A)_{\et}$ obtained from the presheaf $U \mapsto H(U)/J(U)$.

\begin{stmt}
  \label{stmt:quotient-group}
  For each point $t \in \Spec(\A)$, the stalk of the sheaf $H/J$ is
  given by
  \begin{equation*}
    (H/J)_{\bar t} = H(k(\bar{t}))/J(k(\bar{t})).
  \end{equation*}
\end{stmt}

\begin{proof}
  It follows from \cite{milne}*{II Thm. 2.11 and Thm. 2.15} that
  there is an exact sequence
  \begin{equation*}
    1 \to J \to H \to H/J \to 1
  \end{equation*}
  of sheaves of groups on $\Spec(\A)_{\et}$.  Since this sequence is
  exact on stalks, use \ref{stmt:group-scheme-stalk} to see that
  \begin{equation*}
    (H/J)_{\bar{t}}
    \simeq H_{\bar{t}} / J_{\bar{t}} \simeq
    H(\A_{\bar{t}})/J(\A_{\bar{t}}).
  \end{equation*}
  So the assertion will follow once we see that
  $H(\A_{\bar{t}})/J(\A_{\bar{t}}) \simeq
  H(k(\bar{t}))/J(k(\bar{t})).$ Since $H$ is smooth over $S$, the
  natural map $\phi:H(\A_{\bar{t}}) \to H(k(\bar{t}))/J(k(\bar{t}))$
  is surjective \ref{stmt:henselian-lift-to-A}. We have evidently
  $J(\A_{\bar{t}}) \subset \ker \phi$; it remains to see the reverse
  inclusion.

  Recall that $H^o$ is open in $H$. Since $H^o \subset J$, it follows
  that $J \subset H$ is open as well. Since also $J \subset H$ is
  assumed to be closed, one sees that $J$ is a union of connected
  components of $H$.  Let $g \in \ker \phi \subset H(\A_{\bar{t}})$,
  and regard $g$ as a section $\Spec(\A_{\bar{t}}) \to
  H_{/\A_{\bar{t}}}$.  Since $g(\bar{t}) \in J_{\bar{t}}$ and since
  $\Spec(\A_{\bar{t}})$ is connected, it follows that $g \in
  J(\A_{\bar{t}})$. Thus $\ker \phi = J(\A_{\bar{t}})$ as required.
\end{proof}

\subsection{Locally constant component groups}
\label{sub:constant<=>finite-etale}
 Let $X$ be a scheme which is smooth and of finite
type over $\A$. Write $s$ for the closed point of $\Spec(\A)$.
\begin{stmt}
  \label{stmt:lift-to-section}
  If $x \in X(k(\bar{s}))$, choose a finite separable extension $\ell
  \supset k$ such that $x \in X(\ell)$.  Then there is a finite,
  \'etale, local extension $\B\supset \A$ and a section $y \in X(\B)$
  such that $\B$ has residue field $\ell$ and such that $x = y(s') \in
  X(\ell)$, where $s'$ denotes the closed point of $\Spec(\B)$.
\end{stmt}

\begin{proof}
  Since $X$ is of finite type, the existence of the required finite
  separable extension $\ell \supset k$ is immediate.  As in
  \cite{milne}*{I Example 3.4} one may construct a finite, \'etale,
  local extension of $\A$ with residue field
  $\ell$; replacing $\A$ by this extension, we may as well suppose that
  $\ell = k$.

  Now, using \ref{stmt:henselian-lift-to-A} we may find a section $y
  \in X(\Ah)$ over the Henselization $\Ah$ of $\A$ whose image in
  $X(k)$ is $x$.  By construction \cite{milne}*{I \S 4} the
  Henselization $\Ah$ is the limit of \'etale neighborhoods of $\A$;
  the existence of a suitable $\B$ follows at once.
\end{proof}

Recall that $H$ is a group scheme which is smooth, separated, and of
finite type over $\A$.
\begin{stmt}
  \label{stmt:generic-stalk-of-section-enough}
  Let $J \subset H$ be a closed subgroup scheme,
  let $g \in H(\A)$ be a section, and suppose that the image of $g$ in
  $H(K)$ lies in $J(K)$ where $K$ is the field of fractions of $\A$.
  Then $g \in J(\A)$.
\end{stmt}

\begin{proof}
  View $g$ is a morphism $\Spec(\A) \to H$. If $\eta$ denotes the
  generic point of $\Spec(\A)$, the hypothesis means that the
  restriction of $g$ to the dense subset $\{\eta\}$ takes values in
  the closed subset $J$. Since $g$ is continuous, the image of $g$
  must lie in $J$, as required.
\end{proof}

\begin{stmt}
  \label{stmt:injective}
  Let $J \subset H$ be a closed subgroup scheme.  Then the natural map
    \begin{equation*}
      H(\A)/J(\A) \to H(K)/J(K)
    \end{equation*}
    is injective.
\end{stmt}

\begin{proof}
  Suppose that $g \in H(\A)$ and that $gJ(\A)$ is trivial in the group
  $H(K)/J(K)$. Then (the image of) $g$ determines an element of
  $J(K)$.  Since $J$ is closed, it follows from
  \ref{stmt:generic-stalk-of-section-enough} that $g \in J(\A)$.  This
  proves the required injectivity.
\end{proof}


\begin{prop}
  Assume that $J$ is a closed subgroup scheme of $H$ containing $H^o$.
  Then 
  \begin{equation*}
    \#(H/J)_{\bar t} \le \#(H/J)_{\bar \eta} \quad
    \text{for each $t \in \Spec(\A)$},
  \end{equation*}
  where $\eta$ is the generic point of $\Spec(\A)$.  Moreover, the
  following are equivalent:
  \begin{enumerate}
  \item[(i)] The sheaf $H/J$ on $\Spec(\A)_{\et}$ is represented by a
    finite, \'etale $\A$-group scheme.
  \item[(ii)] The sheaf $H/J$ on $\Spec(\A)_{\et}$ is locally constant.
  \item[(iii)] $\#(H/J)_{\bar t}$ is constant on
    $\Spec(\A)$.
  \end{enumerate}
\end{prop}

\begin{proof}
  It suffices to prove the inequality $\#(H/J)_{\bar t} \le
  \#(H/J)_{\bar \eta}$ in case $t$ is the closed point $s \in
  \Spec(\A)$.  Indeed, if $t \in \Spec(\A)$ is arbitrary, one replaces
  $\A$ by the normal local ring $\A_t$; since $t$ is the closed point
  of $\Spec(\A_t)$ one then deduces the required inequality.

  In view of \ref{stmt:quotient-group}, we have $(H/J)_{\bar s} =
  H(k(\bar{s}))/J(k(\bar{s}))$;  choose a finite separable
  extension $\ell \supset k=k(s)$ and elements $x_1,\dots,x_n \in
  H(\ell)$ such that the cosets of the $x_i$ are precisely the
  elements of $H(k(\bar{s}))/J(k(\bar{s}))$.

 We may now use
  \ref{stmt:lift-to-section} to find a finite, \'etale, local
  extension $\B \supset \A$ with residue field $\ell$, and sections
  $y_1,\dots,y_n \in H(\B)$ such that $y_i(s') = x_i$ in $H(\ell)$ for
  $1 \le i \le n$, where $s' \in \Spec(\B)$ is the point lying over $s$.

  Since $J$ is closed in $H$, it follows from \ref{stmt:injective}
  that the classes in $H(k(\bar{\eta}))/J(k(\bar{\eta}))$ of the
  elements $x_i(\eta) \in H(k(\bar{\eta}))$ are all distinct.  Thus
  indeed $\#(H/J)_{\bar{\eta}} \ge n$ as required.

  We now prove that (i) implies (iii); we suppose that there is a
  finite \'etale group scheme $\Gamma$ representing the sheaf $H/J$.
  Since $\A_{\bar x}$ is Henselian, application of \cite{milne}*{I
    Theorem 4.2(c)} shows that the coordinate ring of the finite
  \'etale $\A_{\bar x}$-group scheme $\Gamma_{/\A_{\bar x}}$ is a
  direct product of local rings each with residue field $k(\bar{x})$.
  It follows at once that $\# \Gamma_{\bar\eta} = \# \Gamma_{\bar x}$
  where $\eta$ is the generic point of $\Spec(\A)$, so that indeed
  $\#(H/J)_{\bar x}$ is constant on $\Spec(\A)$.

  We next prove that (iii) implies (ii). Write $n$ for the constant
  value of $\#(H/J)_{\bar x}$. To prove (ii), we must show that $H/J$
  is locally constant. It is enough to prove that $H/J$ is constant
  after we replace $\A$ by a finite, \'etale, local extension; thus,
  we may arrange that there are sections $y_1,\dots,y_n \in H(\A)$ for
  which the cosets $y_i(s)J(k(s))$ are the $n$ distinct elements of
  $H(k(s))/J(k(s)) = (H/J)_{\bar s}$.  It follows from
  \ref{stmt:injective} that the cosets $y_i(\eta)J(K)$ are all
  distinct, where $K = k(\eta)$ is the field of fractions of $\A$.

  Suppose now that $\B \supset \A$ is any finite, \'etale, local
  extension and that $z \in H(\B)$ is any section.  Since by
  assumption $(H/J)_{\bar \eta}$ has $n$ points, we may find $1 \le i
  \le n$ such that $zy_i^{-1} \in J(L)$ where $L$ is the field of
  fractions of $\B$.  But $zy_i^{-1} \in H(\B)$; it then follows from
  \ref{stmt:injective} that $zy_i^{-1} \in J(\B)$. Thus the $y_i$ are
  a full set of coset representatives for the quotient $H(\B)/J(\B)$.

  If now $t \in \Spec(\A)$ is arbitrary, one knows by
  \ref{stmt:lift-to-section} that there is some finite, \'etale, local
  extension $\B \supset \A$ such that $H(\B)/J(\B) \to (H/J)_{\bar t}$
  is surjective. Since the group $(H/J)_{\bar t}$ has exactly $n$
  elements, the above argument shows the natural mapping
  $H(\A)/J(\A) \to (H/J)_{\bar t}$ to be an isomorphism.  It follows that
  $H/J$ is a constant sheaf on $\Spec(\A)_{\et}$, as required.

  The fact that (ii) implies (i) follows from \cite{milne}*{V
    Prop.  1.1} \footnote{Note that \cite{milne}*{Chapter V Prop.
      1.1} is stated for sheaves of Abelian groups, but commutativity
    is not used in the proof.}.
\end{proof}

\subsection{\'Etale central isogenies of group schemes}

Let $H$ and $H_1$ be separated group schemes which are smooth and of
finite type over $\A$ for which $H^o$ and $H_1^o$ are reductive.
Assume that $f:H_1 \to H$ is an \'etale central isogeny over $\A$;
this condition means that the $\A$-morphism $f$ is \'etale, finite,
and faithfully flat, and that $\ker f$ is central in $H$.

\begin{stmt}
  \label{stmt:fo-etale-isog}
  Write $f^o= f_{\mid H_1^o}:H_1^o \to H^o$ for the restriction of $f$
  to $H_1^o$. Then $f^o$ is also an \'etale central isogeny.
\end{stmt}

\begin{proof}
  Since $\ker f$ is central in $H_1$, it is clear that $\ker f^o$ is
  central in $H_1^o$.  Now, $H_1^o$ and $H^o$ are both smooth over
  $\A$. Since the $k(t)$-morphism $f_t:(H_1)_t \to H_t$ is \'etale and
  surjective for each $t \in \Spec(\A)$, it is clear that the same
  holds for $f_t^o:(H_1^o)_t \to {H^o}_t$; thus $f^o$ is \'etale and
  surjective; in particular, $f^o$ is faithfully flat.

  It remains only to show that $f^o$ is finite.  Since $H_1^o$ is
  reductive, we know that $H_1^o \to H_1$ is also a closed immersion
  \ref{stmt:reductive-o-closed}, hence finite. Thus the composition of
  finite morphisms
  \begin{equation*}
    H_1^o \to H_1 \xrightarrow{f} H
  \end{equation*}
  is itself finite; it is then immediate that $f^o$ is finite as well.
\end{proof}

\begin{stmt}
  \label{stmt:H/Ho-iff-H_1/H_1^o}
  The \'etale sheaf $H/H^o$ is represented by a finite \'etale group
  scheme if and only if that is so for $H_1/H_1^o$.
\end{stmt}

\begin{proof}

  One knows that $\ker f$ identifies with the fiber product
  \begin{equation*}
    \begin{CD}
      H_1 @<<< \ker f & = H_1 \times_H \Spec(\A) \\
      @VfVV @VVV \\
      H @<e<< \Spec(\A) \\
    \end{CD}
  \end{equation*}
  where $e$ is the identity section of $H$. Since $f$ is an \'etale
  isogeny, it follows that the subgroup scheme $\ker f \subset H_1$ is
  finite and \'etale over $\Spec(\A)$.

  In view of \ref{stmt:fo-etale-isog}, the same argument shows that
  the subgroup scheme $\ker f^o \subset H_1^o$ is finite and \'etale
  over $\Spec(\A)$ as well.

  Since the sheaves $\ker f$ and $\ker f^o$ are represented by a finite
  \'etale group schemes, one knows by Proposition
  \ref{sub:constant<=>finite-etale} that $\# (\ker f)_{\bar t}$ is
  constant and $\# (\ker f^o)_{\bar t}$ is constant for $t \in
  \Spec(\A)$.  Thus the quotient sheaf $(\ker f)/(\ker f^o)$ has the
  property that $\# ((\ker f)/(\ker f^o))_{\bar t}$ is constant for $t
  \in \Spec(\A)$.

  The sequence of sheaves on $\Spec(\A)_{\et}$
  \begin{equation*}
     1 \to (\ker f)/(\ker f^o) \to H_1/H_1^o \to H/H^o \to 1
  \end{equation*}
  is exact, since it is exact on the stalks of each geometric point
  $\bar{t}$ for $t \in \Spec(\A)$.  It follows at once that
  $\#(H_1/H_1^o)_{\bar t}$ is constant on $\Spec(\A)$ if and only if
  $\#(H/H^o)_{\bar t}$ is constant on $\Spec(\A)$. The result now
  holds by Proposition \ref{sub:constant<=>finite-etale}.
\end{proof}

\subsection{The centralizer of a diagonalizable subgroup scheme}
Let $H$  be as in the previous section; thus
$H$ is smooth and of finite type over $\A$, and $H^o$ is reductive.

We first recall the following result regarding the Weyl group of a
maximal torus of a (connected and) reductive group scheme:
\begin{stmt}{\cite{SGA3}*{Exp. XIX Thm. 2.5}}
  \label{stmt:Weyl-group-finite-and-etale}
  Suppose that $H = H^o$ is reductive with connected geometric fibers,
  let $T \subset H$ be a maximal torus, and let $N_H(T)$ be the
  normalizer in $H$ of $T$. Then the quotient $W = N_H(T)/T =
  N_H(T)/C_H(T)$ is represented by a finite and \'etale group scheme
  over $\A$. In particular, $N_H(T)/T$ is a locally constant sheaf on
  $\Spec(\A)_{\et}$.
\end{stmt}

Suppose now that $D \subset H$ is a closed and smooth \footnote{The
  assumption that $D$ is smooth is imposed here only for lack of
  adequate reference that $L^o = C_H(D)^o$ is reductive; see the
  remark following \ref{stmt:mult-centralizer-in-G}. Given that $L^o$
  is reductive, the proofs of \ref{stmt:L/Lo-when-H-connected} and
  \ref{stmt:L/Lo-via-H/Ho} are independent of the smoothness of $D$.}
subgroup of multiplicative type which is contained in a maximal torus
$T \subset H$, and let $L=C_H(D)$ be the centralizer in $H$ of $D$;
recall \ref{stmt:mult-centralizer-in-G} that $L^o$ is reductive.

\begin{stmt}
  \label{stmt:L/Lo-when-H-connected}
  Assume that $H=H^o$. Then the
  sheaf $L/L^o$ is represented on $\Spec(\A)_{\et}$ by a finite
  \'etale group scheme over $\A$, where $L=C_H(D)$.
\end{stmt}

\begin{proof}
  Since $L^o$ is reductive, one knows by \ref{stmt:reductive-o-closed}
  that $L^o$ is closed in $L$.  If $\eta \in \Spec(\A)$ is the generic
  point, then we have
  \begin{equation*}
    (*)\quad \#(L/L^o)_{\bar \eta} \ge \#(L/L^o)_{\bar t}
  \end{equation*}
  for each $t \in \Spec(\A)$ by the first assertion of Proposition
  \ref{sub:constant<=>finite-etale}. If we show that equality holds in
  $(*)$ for each $t \in \Spec(\A)$, then the desired result follows
  from the equivalence of (i) and (iii) of that same Proposition.

  Recall that $T \subset L^o$ is a maximal torus centralized by $D$.
  Since all maximal tori of $L^o$ are conjugate on the geometric
  fibers ${L^o}_{\bar t}$ for $t \in \Spec(\A)$, the natural map
  \begin{equation*}
    N_L(T)/N_{L^o}(T) \to L/L^o
  \end{equation*}
  determines an isomorphism on each geometric fiber and thus defines an
  isomorphism of sheaves on $\Spec(\A)_{\et}$.

  After replacing $\A$ by a finite, \'etale, local extension, the
  characterization \ref{stmt:Weyl-group-finite-and-etale} shows that
  we may suppose $N_H(T)/T$ to be a constant sheaf. If $n$ denotes the
  (constant) order of the geometric stalks, we may choose sections
  $x_1,x_2,\dots,x_n \in N_H(T)(\A)$ such that $(N_H(T)/T)_{\bar t} =
  \{\overline{x_1(t)},\dots,\overline{x_n(t)}\}$ for each $t \in
  \Spec(\A)$, where $\overline{x_i(t)}$ denotes the coset
  $x_i(t)T(k(\overline{t}))$.

  Now let $y \in N_L(T)(k(\eta))$. Regarding $y$ as an element of
  $N_H(T)(k(\eta))$, we may find $1 \le i \le n$ such that $x_i(\eta)y^{-1} =
  z \in T(K)$.  

  Since the torus $T$ contains $D$, the element $z$ centralizes $D_\eta$;
  it follows that also $x_i(\eta) = yz$ centralizes $D$ -- i.e.
  $x_i(\eta) \in N_L(T)(k(\eta)).$ It now follows from
  \ref{stmt:generic-stalk-of-section-enough} that $x_i \in
  N_L(T)(\A)$.

  We have now proved that the natural map 
  \begin{equation*}
    N_L(T)(\A)/N_{L^o}(T)(\A) \to N_L(T)(k(\eta))/N_{L^o}(T)(k(\eta))
  \end{equation*}
  is surjective. Since for any $t \in \Spec(\A)$ the natural map
  \begin{equation*}
    N_L(T)(\A)/N_{L^o}(T)(\A) \to (N_L(T)/N_{L^o}(T))_{\bar t}
  \end{equation*}
  is injective by the definition of $L^o$, it follows that equality
  holds in $(*)$, as required.
\end{proof}

We will now prove that the assertion of
\ref{stmt:L/Lo-when-H-connected} remains true without the assumption
that $H=H^o$.
\begin{stmt}
  \label{stmt:L/Lo-via-H/Ho}
  If the \'etale sheaf $H/H^o$ is represented by a finite \'etale
  group scheme over $\A$, then $L/L^o$ is represented by
  a finite \'etale group scheme over $\A$ as well.
\end{stmt}

\begin{proof}
  Consider the subgroup $L_1 = C_{H^o}(D)$; we have
  \begin{equation*}
    L^o \subset L_1 \subset L.
  \end{equation*}
  Thus there is an exact sequence of sheaves on $\Spec(\A)_{\et}$
  \begin{equation*}
    1 \to L_1/L^o \to L/L^o \to L/L_1 \to 1.
  \end{equation*}
  It follows from \ref{stmt:L/Lo-when-H-connected} that $L_1/L^o$ is
  represented by a finite \'etale group scheme over $\A$.  
  
  Suppose we show that the order of the geometric fiber $(L/L_1)_{\bar
    t}$ is independent of $t \in \Spec(\A)$.  Using the exactness of
  the above sequence of groups, we see that the order of the geometric
  fiber $(L/L^o)_{\bar t}$ is independent of $t \in \Spec(\A)$.
  Since $L^o$ is reductive and hence closed in $L$ by
  \ref{stmt:reductive-o-closed}, it follows from Proposition
  \ref{sub:constant<=>finite-etale} that $L/L^o$ is locally
  constant and represented by a finite \'etale group scheme over $\A$,
  as required.

  It now remains to prove that $\#(L/L_1)_{\bar t}$ is constant.
  Since $L_1=C_{H^o}(D)$ is closed in $H^o$ and since $H^o$ is closed
  in $H$ by \ref{stmt:reductive-o-closed}, we have that $L_1$ is
  closed in $H$.  Since $L=C_H(D)$ is closed in $H$, $L_1$ is closed
  in $L$ as well. Thus we may apply Proposition
  \ref{sub:constant<=>finite-etale} to study the quotient $L/L_1$.
  That Proposition shows especially that
  \begin{equation*}
    (*)\quad \#(L/L_1)_{\bar \eta} \ge \#(L/L_1)_{\bar t}
  \end{equation*}
  for each $t \in \Spec(\A)$, and the desired result holds if we prove
  that equality holds in $(*)$ for each $t$.

  Recall that $T$ is a maximal torus of $H^o$ containing $D$.  Arguing
  as in \ref{stmt:L/Lo-when-H-connected}, we see that the natural map
  $N_H(T)/N_{H^o}(T) \to H/H^o$ is an isomorphism of sheaves on
  $\Spec(\A)_{\et}$. Since $T \subset L_1$, a similar argument shows
  that the natural map $N_L(T)/N_{L_1}(T) \to L/L_1$ is an isomorphism
  of sheaves on $\Spec(\A)_{\et}$.

  Recall that we have assumed $H/H^o$ to be represented by a finite
  \'etale group scheme; thus $N_H(T)/N_{H^o}(T) \simeq H/H^o$ is
  locally constant.  Since $H^o$ is reductive, also $N_{H^o}(T)/T$ is
  locally constant by \ref{stmt:Weyl-group-finite-and-etale}.  Thus
  after replacing $\A$ by a finite, \'etale, local extension, we may
  suppose that $N_H(T)/N_{H^o}(T)$ and $N_{H^o}(T)/T)$ are constant
  sheaves on $\Spec(\A)_{\et}$. Choose a complete set of
  representatives 
  \begin{equation*}
    x_1,\dots,x_n \in N_H(T)(\A) \quad \text{for the elements of} \quad
    N_H(T)(\A)/N_{H^o}(T)(\A),
  \end{equation*}
  and a complete set of representatives
  \begin{equation*}
      y_1,\dots,y_m \in N_{H^o}(T)(\A) \quad \text{for the elements of} \quad
      N_{H^o}(T)(\A)/T(\A).
  \end{equation*}

  If now $w \in N_L(T)(k(\eta))$, we have $x_i(\eta)y_j(\eta)w^{-1} =
  z \in T(k(\eta))$ for some $1 \le i \le n$ and $1 \le j \le m$. But
  then $x_i(\eta)y_j(\eta) = zy$ centralizes $D_\eta$, hence
  $x_i(\eta)y_j(\eta) \in N_L(T)(k(\eta))$. Since $x_iy_j \in
  N_H(T)(\A)$, it follows from
  \ref{stmt:generic-stalk-of-section-enough} that $x_iy_j \in
  N_L(T)(\A)$.

  This shows  the natural map
  \begin{equation*}
    N_L(T)(\A)/N_{L_1}(T)(\A) \to     N_L(T)(k(\eta))/N_{L_1}(T)(k(\eta))
  \end{equation*}
  to be surjective.

  Notice that the natural map $(L/L_1)_{\bar t} \to (H/H^o)_{\bar t}$
  is injective for each $t \in \Spec(\A)$, hence $L/L_1 \to H/H^o$ is
  an injective mapping of sheaves on $\Spec(\A)_{\et}$. Since the
  natural map $H(\A)/H^o(\A) \to (H/H^o)_{\bar t}$ is injective for
  each $t \in \Spec(\A)$ by the definition of $H^o$, it follows that
  the natural map $L(\A)/L_1(\A) \to (L/L_1)_{\bar t}$ is injective.
  Thus indeed equality holds in $(*)$, as required.
\end{proof}

\section{The component group of a nilpotent centralizer}
\label{section:nilpotent-comp}

Let $\A$ be a local, normal, Noetherian domain, let $G$ be a
$T$-standard reductive group scheme over $\A$. Fix throughout this
section an equidimensional nilpotent section 
\begin{equation*}
  X \in \Lie(G)(\A) =  \glie(\A).
\end{equation*}


Let $C = C_G(X)$ be the centralizer in $G$ of $X$, and recall that $C$
is a smooth group scheme over $\A$; see Proposition
\ref{sub:smooth-over-A}.

\subsection{Replacing $C/C^o$ by $L/L^o$}
According to Theorem \ref{sub:max-tori-and-levi}, $C$ has a Levi
factor locally in the \'etale topology. Thus after replacing $\A$ by a
finite \'etale local extension, we may assume that there is a Levi
factor $L \subset C$.
In this situation, we have the following:
\begin{stmt}
  \label{stmt:C/Co=L/Lo}
  There is an isomorphism $C/C^o \simeq L/L^o$ of sheaves on $\Spec(\A)_{\et}$.
\end{stmt}

\begin{proof}
  It follows from \cite{mcninch-sommers}*{Prop. 12} that the natural
  map defines an isomorphism 
  \begin{equation*}
    (L/L^o)_{\bar t} \simeq (C/C^o)_{\bar
    t} \quad \text{for each $t \in \Spec(\A)$;}
  \end{equation*}
  this means that the natural sheaf map $L/L^o \to C/C^o$ is an
  isomorphism on stalks and is thus an isomorphism.
\end{proof}

\subsection{The adjoint case}
\label{sub:adjoint-case} 
Assume that $G$ is a semisimple group scheme over $\A$, and that $G$ is
of adjoint type -- i.e. that $G_t$ is adjoint for each $t \in
\Spec(\A)$; cf. \S\ref{sub:reductive-group-schemes}.

\begin{theorem}{\cite{mcninch-sommers}}
  The quotient $C/C^o$ is represented on $\Spec(\A)_{\et}$ by a finite
  \'etale group scheme over $\A$. In particular, 
  $(C/C^o)_{\bar t} \simeq (C/C^o)_{\bar{t'}}$ are all $t,t' \in \Spec(\A)$.
\end{theorem}

\begin{proof}
  Replacing $\A$ by a finite, \'etale, local extension we may suppose
  that $C$ has a Levi factor $L$.  According to \ref{stmt:C/Co=L/Lo},
  we know that $C/C^o \simeq L/L^o$.  Since $L^o$ is reductive, $L^o$
  is closed in $L$ \ref{stmt:reductive-o-closed}. Thus we may apply
  Proposition \ref{sub:constant<=>finite-etale}. According to that
  Proposition, the Theorem will follow once we know that
  $|(L/L^o)_{\bar t}|$ is constant for $t \in \Spec(\A)$.  Since for
  each $t \in \Spec(\A)$, we have assumed $G_{\bar t}$ to be
  semisimple and adjoint, that constancy follows from
  \cite{mcninch-sommers}*{Theorem 36}.
\end{proof}

\subsection{\'Etale central isogenies}
\label{sub:component-etale-isog}
Let $G$ and $G_1$ be $T$-standard reductive group schemes over
$\A$, and let $f:G \to G_1$ be an \'etale central isogeny. This means
that $f$ is finite, \'etale, and faithfully flat, and $\ker f$ is
central in $G$.

\begin{stmt}
  $df:\Lie(G)(\A) \to \Lie(G_1)(\A)$ is an isomorphism of $\A$-modules.
\end{stmt}

\begin{proof}
  Since $f_t:G_t \to G_{1,t}$ is a separable central isogeny for each
  $t \in \Spec(\A)$, one knows that $\dim G_t = \dim G_{1,t}.$ Thus,
  $\Lie(G)(\A)$ and $\Lie(G_1)(\A)$ are free $\A$-modules of the same
  rank. Writing $\mathfrak{m}$ for the unique maximal ideal of $\A$,
  one knows that
  \begin{equation*}
    \Lie(G)(k(s)) = \Lie(G)(\A)/\mathfrak{m}\Lie(G)(\A),
  \end{equation*}
  with a similar statement for $G_1$.  Since $df_s:\Lie(G)_s \to
  \Lie(G_1)_s$ is an isomorphism, it follows from the Nakayama lemma
  that $df:\Lie(G)(\A) \to \Lie(G_1)(\A)$ is an isomorphism.
\end{proof}

Let $X \in \Lie(G)(\A)$ and $X_1 \in \Lie(G_1)(\A)$ be nilpotent
sections, and suppose that $df(X) = df(X_1)$.  Write $C = C_G(X)$, $N
= N_G(X)$, $C_1 = C_{G_1}(X_1)$ and $N_1 = N_{G_1}(X_1)$. 

\begin{stmt}
  $X$ is an equidimensional nilpotent section of $\Lie(G)$ if and only
  if $X_1$ is an equidimensional nilpotent section of $\Lie(G_1)$.
\end{stmt}

\begin{proof}
  Indeed, it is clear for each $t$ that $f_t$ restricts to a
  separable isogeny
  \begin{equation*}
    f_{t\mid C_t}:C_t \to C_{1,t}
  \end{equation*}
  of $k(t)$-group schemes, whence $\dim C_t = \dim C_{1,t}$.
\end{proof}

\begin{stmt}
  \label{stmt:L-->L_1}
  Locally in the \'etale topology there are Levi factors $L
  \subset C$ and $L_1 \subset C_1$ for which $f_{\mid L}$ determines a
  finite, \'etale, and faithfully flat map of group schemes $f_{\mid
    L}:L \to L_1$.
\end{stmt}

\begin{proof}
  After possibly replacing $\A$ by a finite, \'etale, local extension,
  we may use Theorem \ref{sub:max-tori-and-levi}(a) to find a
  homomorphism $\phi:\Gm \to G$ such that $\phi_t$ is a cocharacter of
  $G_t$ associated to $X(t)$ for each $t \in \Spec(\A)$. If $\psi = f
  \circ \phi$, it follows from \cite{mcninch-rat}*{Lem. 14} that
  $\psi_t$ is a cocharacter of $G_{1,t}$ associated with $X_1(t)$ for
  each $t$,

  Using \ref{stmt:phi=>Levi}, one knows that the centralizer $L$ of
  the image of $\phi$ in $C$ is a Levi factor, and the centralizer
  $L_1$ of the image of $\psi$ in $C_1$ is a Levi factor.

  It is clear that $f$ restricts to a morphism $f_{\mid L}:L \to L_1$;
  we only must argue that $f_{\mid L}$ is finite, \'etale, and
  faithfully flat. 

  For that, we notice first that $f_{t\mid L_t}:L_t \to L_{1,t}$ is a
  separable $k(t)$-isogeny for each $t \in \Spec(\A)$. It follows at
  once that $f$ is faithfully flat. Moreover, since $L$ and $L_1$ are
  smooth over $\A$, \cite{SGA1}*{Exp. II, Cor. 2.2} shows that $f_{\mid
    L}$ is smooth. Since $f$ is finite, it follows that $f_{\mid L}$
  is quasi-finite. But then $f_{\mid L}$ is \'etale 
  \cite{SGA1}*{Exp. II, Cor. 1.4}.

  It remains to show that $f_{\mid L}$ is a finite morphism.  Note
  first that the inclusions $C \subset G$ and $C_1 \subset G_1$ are
  closed embeddings (see \S\ref{sub:smoothness}), and the inclusions
  $L \subset C$ and $L_1 \subset C_1$ are closed embeddings
  \ref{stmt:mult-centralizer}.  Since $f$ is finite hence proper
  \cite{milne}*{I Prop. 1.4}, the composition $L \subset C \subset G
  \xrightarrow{f} G_1$ is a proper map. Since $L_1 \subset G_1$ is a
  closed embedding, it follows from \cite{liu}*{Prop.  3.3.16} that
  $f_{\mid L}:L \to L_1$ is proper.  Since $f_{\mid L}$ is quasifinite
  and proper, $f_{\mid L}$ is finite by \cite{milne}*{I Cor. 1.1}.
\end{proof}

\begin{prop}
  The sheaf $C/C^o$ is represented on $\Spec(\A)_{\et}$ by a finite
  \'etale $\A$-scheme if and only if that is so for $C_1/C^o_1$.
\end{prop}

\begin{proof}
  Replacing $\A$ be a finite, \'etale local extension, we may use
  \ref{stmt:L-->L_1} to find Levi factors $L \subset C$ and $L_1
  \subset C_1$ for which $f_{\mid L}$ determines a finite, \'etale,
  faithfully flat morphism $f:L \to L_1$. Then we have
  \begin{equation*}
    C/C^o \simeq L/L^o \quad \text{and} \quad C_1/C^o_1 \simeq L_1/L_1^o
  \end{equation*}
  by \ref{stmt:C/Co=L/Lo}. Since $L^o$ and $L_1^o$ are reductive, the
  Proposition now follows by applying \ref{stmt:H/Ho-iff-H_1/H_1^o}.
\end{proof}

\subsection{The centralizer of a torus}
\label{sub:component-torus-centralizer}

Let $S \subset G$ be a torus, and let $M = C_G(S)$ be the centralizer
of $S$ in $G$. Then $M$ is a reductive group scheme over $\A$ with
connected geometric fibers.  Suppose that $X \in \Lie(M)(\A)$ is an
equidimensional nilpotent section, and  write $C = C_G(X)$ and $C_M =
C_M(X)$.

\begin{stmt}
  \label{stmt:equidim-for-G}
  $X$ is an equidimensional nilpotent section of $\Lie(G)(\A)$ as
  well. In particular,  $C$ is equidimensional and hence smooth.
\end{stmt}

\begin{proof}
  It follows from Corollary \ref{sub:max-tori-and-levi} that there are
  subgroup schemes $Q \subset L \subset M$ such that $L$ is a Levi
  subgroup scheme of a parabolic subgroup scheme of $M$, $Q$ is a
  parabolic subgroup scheme of $M$, and $(L_{\bar t},Q_{\bar t})$ is
  the Bala-Carter datum of $X(t)$ for each $t \in \Spec(\A)$.

  Now use Theorem \ref{sub:existence-equidim} to find an
  equidimensional nilpotent section $Y \in \Lie(G)(\A)$ whose
  Bala-Carter datum coincides with that of $X(\eta)$. Since the
  Bala-Carter datum of $X(t)$ in $\Lie(G)_t$ is determined by the
  Bala-Carter datum of $X(t)$ in $\Lie(M)_t$, it follows from the
  Bala-Carter theorem that $X(t)$ and $Y(t)$ are conjugate by an
  element of $G(k(\bar{t}))$.  Since $Y$ is equidimensional, it
  follows that $X$ is equidimensional as well.
\end{proof}

\begin{prop}
    If the sheaf $C/C^o$ is represented on $\Spec(\A)_{\et}$ by a
    finite \'etale $\A$-group scheme, the same holds for $C_M/C^o_M$.
\end{prop}

\begin{proof}


  In view of \ref{stmt:equidim-for-G}, one knows that $C$ and $C_M$
  are smooth over $\A$. Thus the centralizers $L \subset C$ and $L_1
  \subset C_M$ of the image of $\phi$ are closed subgroup schemes
  which are smooth over $\A$.  Using \ref{stmt:phi=>Levi}, one knows
  that $L$ is a Levi factor in $C$ and that $L_1$ is a Levi factor in
  $C_M$. Now,
  \begin{equation*}
    C/C^o \simeq L/L^o \quad \text{and} \quad C_M/C_M^o \simeq L_1/L_1^o
  \end{equation*}
  by \ref{stmt:C/Co=L/Lo}. 

  Since $\phi$ evidently centralizes the torus $S$, it is clear
  that $S \subset L_1^o \subset L^o$.  Since the centralizer of $S$ in
  $L^o$ is a reductive subgroup scheme, we may find a maximal torus $T
  \subset L^o$ containing $S$ -- use
  \ref{stmt:reductive-has-max-torus} and \ref{stmt:D-contained-in-T}.
  The Proposition now  follows from \ref{stmt:L/Lo-via-H/Ho}.
\end{proof}

\begin{rem}
  It is not clear -- to the author, at least -- whether the
  Proposition is true when $T$ is replaced by any diagonalizable
  subgroup scheme $D \subset G$; in the notation of the (proof of the)
  Proposition, the difficulty lies in the fact that $D$ need not be
  contained in a maximal torus of $L_1^o$, so that
  \ref{stmt:L/Lo-via-H/Ho} is inadequate.
\end{rem}

\subsection{The component group of $C$}
\label{sub:main-result-for-component-group}
Let $G$ be a $T$-standard reductive group scheme over $\A$.  Let $X
\in \glie(\A)$ be an equidimensional nilpotent section, let $C =
C_G(X)$, and assume that the pair $(G,X)$ is allowable.
\begin{theorem}
  The \'etale sheaf $C/C^o$ is represented on $\Spec(\A)_{\et}$ by a
  finite \'etale group scheme over $\A$. In particular, 
  $(C/C^o)_{\bar{t}} \simeq (C/C^o)_{\bar{t'}}$ for all $t,t' \in \Spec(\A)$.
\end{theorem}

\begin{proof}
  Let first $G$ be semisimple and assume the fiber characteristics are
  all very good for $G$. After replacing $\A$ by a finite, \'etale,
  local extension, we may suppose that $G$ is split; let $G_{\ad}$ be
  the corresponding group of adjoint type and let $f:G \to G_{\ad}$ be
  the corresponding map \ref{stmt:adjoint-group}.  In view of our
  assumptions, \ref{stmt:very-good-p-adjoint-group} shows that $f$ is
  an \'etale central isogeny. Since the assertion of the Theorem holds
  for the pair $G_{\ad},df(X)$ by Theorem \ref{sub:adjoint-case},
  the assertion for $G$ now follows from Proposition
  \ref{sub:component-etale-isog}.

  It is then clear that the assertion of the Theorem holds when $G$ is
  a group of the form $H=H_1 \times S$ where $S$ is a torus and where
  $H_1$ is semisimple and the characteristic of $k(t)$ is very good
  for $H_t$ for each $t \in \Spec(\A)$.  If $S_0$ is a torus in the
  group $H$, Proposition \ref{sub:component-torus-centralizer} shows
  that the Theorem holds for $M = C_H(S_0)$. 

  If $G$ is any $T$-standard group, there is an \'etale isogeny
  between $G$ and a group of the form $M$ as above; thus the assertion
  of the Theorem follows from another application of Proposition
  \ref{sub:component-etale-isog}.
\end{proof}

\subsection{Proof of Theorem \ref{theorem:component}}
\label{sub:proof-component}

Recall that $E$ is an algebraically closed field of characteristic 0,
and that $F$ is an algebraically closed field of characteristic $p>0$.
Theorem \ref{theorem:component} is a consequence of the following
more general result:
\begin{theorem}
  Let $G_E$ and $G_F$ be reductive groups respectively over $E$ and
  over $F$, assume that the root datum of $G_E$ coincides with that of
  $G_F$, and assume that $G_F$ is $T$-standard. Let $X_E \in \glie_E$,
  $X_F \in \glie_F$ be nilpotent elements with the same Bala-Carter
  data, and let $C_E$ and $C_F$ be their respective centralizers.  Then
  $C_E/C_E^o \simeq C_F/C_F^o$.
\end{theorem}

\begin{proof}
  As in the proof of Theorem
  \ref{sub:proof-levi-factor}, let $\A$ be the ring of Witt vectors
  \footnote{The same remarks concerning the choice of $\A$ made in the
    footnote in \S\ref{sub:proof-levi-factor} apply here as well; we
    could instead take for $\A$ any normal, local, Noetherian domain
    with infinite residue field of characteristic $p$ and field of
    fractions of characteristic $0$.}  \cite{serre-local-fields}*{II
    \S6} with residue field an algebraic closure of the finite field
  $\mathbf{F}_p$.  Using \ref{stmt:choose-favorite-fields}, we see --
  as in the proof of Theorem \ref{sub:proof-levi-factor} -- that it
  suffices to prove the Theorem after replacing
  $F$ by the residue field of $\A$ and $E$ by an algebraic closure of
  the field of fractions of $\A$, and after replacing $X_F$ and $X_E$
  by nilpotent elements with the given Bala-Carter datum.

  Again, let $G$ be a split reductive group scheme over $\A$ with the
  given root datum.  Using Theorem \ref{sub:existence-equidim} and the
  Bala-Carter Theorem, we may suppose that there is an equidimensional
  nilpotent section $X$ for which $X_F$ is conjugate to $X(s)$ and for
  which $X_E$ is geometrically conjugate to $X(\eta)$.

  If $C = C_G(X)$ denotes the centralizer subgroup scheme, it now follows from
  Theorem \ref{sub:main-result-for-component-group} that
  $(C/C^o)_{\bar s} \simeq (C/C^o)_{\bar \eta}$. Thus, the component
  group of the centralizer in $G_{\bar \eta}$ of $X(\eta)$ is
  isomorphic to the component group of the centralizer in $G_{\bar s}$
  of $X(s)$, as required.
\end{proof}


\begin{bibdiv}

  \begin{biblist}[\resetbiblist{XXXXXX}]

  

  \bib{alek}{article}{
    author = {A.~V. Alekseevski{\u\i}},
    title = {Component
      groups of centralizers of unipotent elements in semisimple
      algebraic groups},
    journal = {Akad. Nauk Gruzin. SSR Trudy Tbiliss. Mat.
      Inst. Razmadze},
    volume = {62},
    year = {1979},
    note = {Collection of articles on
      algebra, 2.},
    label = {Al 79}
    }

  \bib{borel-LAG}{book}{
      author = {Borel, Armand},
      title = {Linear Algebraic Groups},
      year = {1991},
      publisher = {Springer Verlag},
      volume = {126},
      series = {Grad. Texts in Math.},
      edition = {2nd ed.},
      label = {Bor 91}
    }



    \bib{carter}{book}{
      author = {Roger~W. Carter},
      title =  {Finite groups of {L}ie type: conjugacy classes and
        complex characters},
      publisher = {John Wiley \& Sons Ltd.},
      place = {Chichester},
      year = {1993},
      note = {Reprint of the 1985 original},
      label = {Ca 93}
    }

    \bib{EGA-IV}{article}{
      author={Grothendieck, A.},
      title={\'El\'ements de g\'eom\'etrie alg\'ebrique. IV. \'Etude locale
        des sch\'emas et des morphismes de sch\'emas. III},
   journal={Inst. Hautes \'Etudes Sci. Publ. Math.},
    number={28},
      date={1966},
     pages={255},
      issn={0073-8301},
      label={EGA IV}}

    \bib{SGA1}{collection}{ 
      title={Rev\^etements \'etales et groupe  fondamental (SGA 1)}, 
      series={Documents Math\'ematiques 3},
      author={A. Grothendieck}, 
      note={S\'eminaire de G\'eometrie  Alg\'ebrique du Bois Marie. [Updated and annotated reprint of
        the 1971 original publication in Lecture Notes in Math., 224, Springer,
        Berlin.]}, 
      publisher={Soci\'et\'e Math\'ematique de France},
      place={Paris}, 
      date={2003}, 
      label={SGA 1} }


      \bib{SGA3}{book}{
      author = {Grothendieck, A.},
      author = {Demazure, M.},
      title  = {Sch\'emas en Groupes (SGA 3). I, II, III},
      series   = {Lectures Notes in Math.},
      volume =  {151, 152, 153},
      publisher = {Springer Verlag},
      note = {S\'eminaire de G\'eometrie Alg\'ebrique du Bois Marie},
      place  = {Heidelberg},
      year  = {1965},
      label = {SGA 3}    }


  \bib{Hum95}{book}{ 
    author={Humphreys, James~E.},
    title={Conjugacy classes in semisimple algebraic groups},
    series={Math. Surveys and Monographs}, 
    publisher={Amer. Math. Soc.}, 
    date={1995}, 
    volume={43},
    label={Hu 95}}


    \bib{jantzen-nil}{article}{ 
      author={Jantzen, Jens~Carsten},
      book = {
        title = {Lie Theory: Lie Algebras and Representations},
        series = {Progress in Mathematics},
        publisher = {Birkh\"auser},
        editor = {Anker, J-P},
        editor = {Orsted, B},
        place = {Boston},
        volume = {228},
        date = {2004}      },
      title={Nilpotent orbits in representation theory}, 
      pages = {1\ndash211},
      label = {Ja 04}}

    

    \bib{JRAG}{book}{
      author={Jantzen, Jens Carsten},
      title={Representations of algebraic groups},
      series={Mathematical Surveys and Monographs},
      volume={107},
      edition={2},
      publisher={American Mathematical Society},
      place={Providence, RI},
      date={2003},
      pages={xiv+576},
      label = {Ja 03}
    }
	

  \bib{kempf-instab}{article}{
    author={Kempf, George~R.},
     title={Instability in invariant theory},
      date={1978},
      ISSN={0003-486X},
   journal={Ann. of Math. (2)},
    volume={108},
    number={2},
     pages={299\ndash 316},
    label ={Ke 78} }      

    \bib{KMRT}{book}{ 
      author={Knus, Max-Albert}, 
      author={Merkurjev,  Alexander}, 
      author={Rost, Markus}, 
      author={Tignol, Jean-Pierre},
      title={The book of involutions}, 
      series={Amer. Math. Soc. Colloq.  Publ.}, 
      publisher={Amer. Math. Soc.}, 
      date={1998}, 
      volume={44},
      label = {KMRT}  }



    \bib{liu}{book}{
      author = {Liu, Qing},
      title = {Algebraic geometry and arithmetic curves},
      note = {Translated from the French by Reinie Ern{\'e}},
      series={Oxford Graduate Texts in Mathematics},
      number = {6},
      publisher={Oxford University Press},
      year = {2002},
      label = {Li 02}}


    

    \bib{mcninch-sommers}{article}{
      author={McNinch, George J.},
      author={Sommers, Eric},
      title={Component groups of unipotent centralizers in good
        characteristic},
      note={Special issue celebrating the 80th birthday of Robert
        Steinberg. math.RT/0204275},
      journal={J. Algebra},
      volume={260},
      date={2003},
      number={1},
      pages={323\ndash 337},
      issn={0021-8693},
      label = {MS 03}}

 \bib{mcninch-testerman}{article}{
      author={McNinch, George~J.},
      author={Testerman, Donna~M.},
      title={Completely reducible $\operatorname{SL}(2)$-homomorphisms},
      journal = {Transact. AMS},
      volume = {359},
      number = {9},
      year = {2007},
      pages = {4489\ndash 4510}
      label={MT 07}}

    \bib{mcninch-rat}{article}{
      author = {McNinch, George~J.},
      title = {Nilpotent orbits over ground fields of good characteristic},
      volume = {329},
      pages = {49\ndash 85},
      year = {2004},
      note={arXiv:math.RT/0209151},
      journal = {Math. Annalen},
      label={Mc 04}}

    \bib{mcninch-optimal}{article}{
      author={McNinch, George~J.},
      title = {Optimal $\operatorname{SL}(2)$-homomorphisms},
      date = {2005},
      volume = {80},
      pages = {391 \ndash 426},
      journal = {Comment. Math. Helv.},
      note = {arXiv:math.RT/0309385},
      label = {Mc 05}}

    \bib{mcninch-nil-sum}{article}{
      author = {McNinch, George~J.},
      title = {On the centralizer of the sum of commuting nilpotent elements},
      date = {2006},
      journal = {J. Pure and Applied Algebra},
      volume = {206},
      pages = {123\ndash140},
      year = {2006},
      note = {arXiv:math.RT/0412283},
      label = {Mc 06}}

     \bib{milne}{book}{
       author = {Milne, James},
       title = {{\'E}tale Cohomology},
       year = {1980},
       publisher = {Princeton University Press},
       label={Mil 80}    }


     \bib{mizuno}{article}{
       author={Mizuno, Kenzo},
       title = {The conjugate classes of unipotent elements of the
         {C}hevalley groups ${E}\sb{7}$\ and ${E}\sb{8}$},
       year = {1980},
       journal = {Tokyo J. Math.},
       volume = {3},
       number = {2},
       pages = {391\ndash461},
       label={Miz 80}
       }


     \bib{premet}{article}{
       author={Premet, Alexander},
       title={Nilpotent orbits in good characteristic and the Kempf-Rousseau theory},
       note={Special issue celebrating the 80th birthday of Robert
         Steinberg.},
       journal={J. Algebra},
       volume={260},
       date={2003},
       number={1},
       pages={338\ndash 366},
       label = {Pr 03}}

     \bib{richardson}{article}{
       author={Richardson, R. W.},
       title={On orbits of algebraic groups and Lie groups},
       journal={Bull. Austral. Math. Soc.},
       volume={25},
       date={1982},
       number={1},
       pages={1--28},
       label = {Ri 82}
     }

     \bib{serre-local-fields}{book}{
       author = {Jean-Pierre Serre},
       title = {Local Fields},
       year = {1979},
       publisher = {Springer Verlag},
       label = {Se 79}}

    \bib{serre-sem-bourb}{article}{
      author = {Serre, Jean-Pierre},
      title =  {Compl{\`e}te R{\'e}ductibilit{\'e}},
      publisher={Soci\'et\'e Math\'ematique de France},
      note = {S{\'e}minaire Bourbaki 2003/2004},
      pages= {Expos{\'e}s  924-937, pp. 195\ndash217},
      journal = {Ast{\'e}risque},
      volume = {299},
      year = {2005},
      label={Ser 05}      }



    \bib{sommers}{article}{
      author = {Sommers, Eric},
      title = {A generalization of the {B}ala-{C}arter theorem for
        nilpotent orbits},
      journal = {Internat. Math. Res. Notices},
      year = {1998},
      pages = {539\ndash562},
      number = {11},
      label = {So 98}
      }

    \bib{springer-LAG}{book}{ 
      author={Springer, Tonny~A.}, 
      title={Linear algebraic groups}, 
      edition={2}, 
      series={Progr. in Math.},
      publisher={Birkh{\"a}user}, address={Boston}, date={1998},
      volume={9}, 
      label={Sp 98}}

    \bib{steinberg-endomorphisms}{book}{
      author={Steinberg, Robert},
      title={Endomorphisms of linear algebraic groups},
      series={Memoirs of the American Mathematical Society, No. 80},
      publisher={American Mathematical Society},
      place={Providence, R.I.},
      date={1968},
      pages={108},
      label = {St 68}} 
  \end{biblist}
\end{bibdiv}
\end{document}